\documentclass[11pt]{amsart}

\usepackage{amsfonts, enumerate} 
\usepackage{amssymb}
\usepackage[final]{graphicx, graphics,stmaryrd} 

\usepackage[margin=3cm]{geometry} 
\usepackage{amscd}
\usepackage[all]{xy}

\usepackage[breaklinks,bookmarksopen,bookmarksnumbered]{hyperref}

\usepackage{tensor}
\usepackage[english]{babel}
\usepackage[utf8]{inputenc}
\usepackage{amsthm}
\usepackage{graphics}
\usepackage{amsmath}
\usepackage{amstext}
\usepackage{engrec}
\usepackage{rotating}
\usepackage[safe,extra]{tipa}
\usepackage{multirow}
\usepackage{enumitem}
\usepackage{bm}
\usepackage{stackrel} 
\usepackage{tikz-cd} 
\numberwithin{equation}{section}

\newtheorem{teo}{Theorem}[subsection]

\newtheorem{prop}[teo]{Proposition}
\newtheorem{lemma}[teo]{Lemma}

\newtheorem*{teon}{Theorem}

\theoremstyle{definition}
\newtheorem{defin}[teo]{Definition}
\newtheorem{rmk}[teo]{Remark}

\newcommand{\mb}{\mathbb}
\newcommand{\mc}{\mathcal}
\newcommand{\msf}{\mathsf}
\newcommand{\mfk}{\mathfrak}

\def\Oo{\mathcal O}
\def\R{\mathcal{R}}
\def\C{\mathbb C}

\def\Z{\mathbb Z}

\def\Q{\mathbb Q}
\def\p{\mathbb P}

\def\E{\mathbb E}
\def\M{\mathcal{M}}
\def\A{\mathbb A}

\def\mm{\bm{\mathsf{M}}}
\def\rr{\bm{\mathsf{R}}}

\def\H{\mfk{H}}

\newcommand{\ssC}{\mathsf{C}}
\newcommand{\ssL}{\mathsf{L}}

\newcommand{\m}{\mbox}

\newcommand{\cor}{\textit}

\newcommand{\fine}{\qed\newline}
\newcommand{\xx}{\otimes}

\newcommand{\beginCD}{\begin{equation*}\begin{CD}}
\newcommand{\enCD}{\end{CD}\end{equation*}}

\newcommand{\G}{\mfk G}

\newcommand{\we}{\wedge}

\DeclareMathOperator{\age}{age}

\DeclareMathOperator{\defo}{Def}
\DeclareMathOperator{\Aut}{Aut}
\DeclareMathOperator{\diag}{Diag}

\DeclareMathOperator{\ord}{ord}
\DeclareMathOperator{\im}{Im}
\DeclareMathOperator{\id}{id}
\DeclareMathOperator{\sing}{Sing}

\DeclareMathOperator{\Ho}{Hom}
\DeclareMathOperator{\Ker}{Ker}
\DeclareMathOperator{\Spec}{Spec}

\DeclareMathOperator{\GL}{GL}

\DeclareMathOperator{\QR}{QR}

\DeclareMathOperator{\nc}{nc}

\def\rr{\bm{\mathsf{R}}}

\def\lq{\llbracket}
\def\rq{\rrbracket}

\DeclareMathOperator{\cov}{Cov}
\newcommand{\nor}{\msf {nor}}

\DeclareMathOperator{\cord}{c-ord}

\DeclareMathOperator{\bal}{bal}

\DeclareMathOperator{\adm}{Adm}
\def\E{\mb E}

\def\rr{\bm{\mathsf{R}}}
\def\mmu{\bm{\mu}}

\DeclareMathOperator{\gen}{gen}
\DeclareMathOperator{\sep}{sep}

\DeclareMathOperator{\lcm}{lcm}
\DeclareMathOperator{\Def}{Def}
\DeclareMathOperator{\Fun}{Fun}

\def\t{\tau}

\DeclareMathOperator{\inn}{Inn}

\DeclareMathOperator{\sub}{Sub}
\DeclareMathOperator{\Zar}{Zar}

\title{Moduli of $G$-covers of curves: geometry and singularities}
\author{Mattia, Galeotti}
\address{Università di Bologna,\ Piazza di Porta S. Donato, 5, 40126 Bologna, Italy}
\email{galeotti.mattia.work@gmail.com}

\begin{document}
\maketitle

\begin{abstract}
We analyze the singular locus and the locus of non-canonical singularities
of the moduli space $\overline\R_{g,G}$ of curves with a $G$-cover for any finite group $G$.
We show that non-canonical singularities are of two types:
$T$-curves, that is singularities lifted from the moduli space $\overline\M_g$ of stable curves,
and $J$-curves, that is new singularities entirely characterized by the dual graph of the cover.
Finally, we prove that in the case $G=S_3$, the $J$-locus is empty, which is the first
fundamental step in evaluating the Kodaira dimension of~$\overline\R_{g,S_3}$.
\end{abstract}

\section{Introduction}
This is the first of two papers whose goal is to analyze the birational geometry
of the moduli space of curves equipped with a $G$-cover, where $G$ is any finite group.
In particular we focus on the case $G=S_3$, the symmetric group of order $3$. 

The moduli space $\M_g$ of smooth curves
of genus $g$ is a widely studied object along with its Deligne-Mumford compactification $\overline\M_g$
described for the first time in \cite{delmum69}.
This compactification is the moduli space
of genus $g$ stable curves, that is curves admitting nodal singularities  and a finite
automorphism group. The birational geometry of $\overline\M_g$ was first approached Eisenbud, Harris and Mumford~\cite{harmum82, har84, eihar87},
proving that it is a variety of general type for genus $g> 23$.
The cases $g=22,\ 23$ were recently solved by Farkas-Jensen-Payne \cite{fjp20}, 
proving that $\overline\M_{22}$ and $\overline\M_{23}$ are of general type too.

The present work fits in the framework of
finite covers of $\overline\M_g$,
whose study is motivated by the fact that in many
cases the transition to the general type happens for genus lower than~$22$.
Farkas and Verra (see \cite{farver10}) focused in the case of odd spin curves; Chiodo-Eisenbud-Farkas-Schreyer work \cite{cefs13}
analyzes the moduli of curves with a $3$-torsion bundles;
in both cases the moduli space is of general type for $g\geq 12$.
For this type of results is fundamental an analysis of the singular locus.
This has been done by Chiodo and Farkas in \cite{chiofar12} for curves with an $\ell$-torsion bundle,
also called level $\ell$ curves. In his work \cite{gale15}, the author generalized this analysis to the case
of the moduli space $\R_{g,\ell}^k$ of curves with a line bundle $L$ such that $L^{\xx\ell}\cong\omega^{\xx k}$.

Here we propose another generalization of Chiodo and Farkas approach, by treating curves with a $G$-cover
for any finite group $G$, where
the case of level $\ell$ curves is equivalent to $G=\mmu_{\ell}$ a cyclic group.
In order to compactify the moduli space $\R_{g,G}$ of genus $g$ smooth curves with a principal $G$-bundle, 
we introduce 
two notions of covers: twisted $G$-covers and admissible $G$-covers.
Twisted covers are treated in \cite{chiofar12} as balanced maps $\phi\colon\msf C\to BG$ where $\msf C$ is a twisted curve,
that is a Deligne-Mumford stack whose coarse space is a stable curve and with non-trivial cyclic stabilizer
at some nodes. For a wide introduction to twisted curves and their moduli see for example \cite{acv03, abravis02, chio08}.
Admissible $G$-covers are principal $G$-bundles admitting ramification points over some nodes. The two cover notions
are proved equivalent in \cite{acv03}, as recalled here in Theorem \ref{teo_3e}.

The main result we propose is the description of the singular locus $\sing\overline\R_{g,G}$
and the non-canonical singular locus $\sing^{\nc}\overline\R_{g,G}$. In particular, we are interested
in characterizing the singularities outside the preimage of singular points of $\overline\M_g$.
In order to achieve this, for any twisted $G$-cover $(\ssC,\phi)$ we consider the group $\underline\Aut_C(\ssC,\phi)$
of ghost automorphisms, \cor{i.e.}~$\ssC$ automorphisms lifting to $\phi$ and acting trivially on the
coarse curve $C$. As any singularity of $\overline\R_{g,G}$ is a quotient singularity, there
are some tools allowing its description, such as quasireflections (see Definition \ref{def_qr})
and the age invariant, in particular via the notion of \cor{junior} group (see Definition \ref{def_junior}).
We denote by $\QR\subset\underline\Aut_C(\ssC,\phi)$ the subgroup generated by quasireflections.
Moreover, if $\pi\colon \overline\R_{g,G}\to \overline\M_g$ is the natural projection,
we denote by $N_{g,G}:=\sing\overline\R_{g,G}\cap\pi^{-1}\sing\overline\M_g$ and
$T_{g,G}:=\sing^{\nc}\overline\R_{g,G}\cap\pi^{-1}\sing^{\nc}\overline\M_g$ the loci of singularities
lifted from $\overline\M_g$. Theorems \ref{teo_smooth} and \ref{teo_jt}, summarized below, 
say that the ``new'' singularities are characterized by their ghost structure.

\begin{teon} 
If $H_{g,G}\subset \overline\R_{g,G}$ is the locus of twisted $G$-covers $(\ssC,\phi)$ such
that $\underline\Aut_C(\ssC,\phi)$ is not generated by quasireflections, then
$$\sing\overline\R_{g,G}=N_{g,G}\cup H_{g,G}.$$
If $J_{g,G}\subset\overline\R_{g,G}$ is the locus such that $\underline\Aut_C(\ssC,\phi)\slash\QR$ is
a junior group, then
$$\sing^{\nc}\overline\R_{g,G}=T_{g,G}\cup J_{g,G}.$$
\end{teon}

In order to approach the problem of evaluating the Kodaira dimension of $\overline\R_{g,G}$,
a fundamental step is proving the pluricanonical form extension result, similarly to what has been done
for $\overline\M_g$ in \cite{harmum82}. As last result we prove in Theorem \ref{teo_ncs3}
that the $J$-locus is empty for $G=S_3$, and this will be the starting point to the extension
of pluricanonical forms over a desingularization of $\overline\R_{g,S_3}$,
because it allows the generalization of Harris-Mumford techniques.
\begin{teon}
In the case of the symmetric group $G=S_3$, the $J$-locus $J_{g,S_3}$ is empty
for any genus $g\geq 2$. Therefore $\sing^{\nc}\overline\R_{g,S_3}=T_{g,S_3}$.
\end{teon}

As a direct application, in our next paper we are going to prove that the moduli
space of genus $g$ connected twisted $S_3$-covers is of general type
for any odd genus $g\geq 13$.

In section \S\ref{moduli1} we introduce  the different notions of covers and recall their equivalence.
In \S\ref{chap_mod} we review the dual graph and torsor notions, they are very important in describing
the structure of twisted covers and their ghost automorphisms. In \S\ref{soR} and \S\ref{ncs}
we prove the main results for the loci of singularities.

\section*{Acknowledgements}
This work is the completion of my PhD thesis at \cor{Institut de Mathématiques de Jussieu}, thus  I want to thank
my PhD supervisor Alessandro Chiodo, for his patience and his advices. Furthermore, I want to thank
Michele Bolognesi and Gavril Farkas for the careful reading and the comments on my thesis work. Finally, I am grateful
to Claudio Fontanari, my post-doc tutor at \cor{Università degli Studi di Trento}, for his important suggestions,
and to Roberto Pignatelli for the fruitful conversations.\newline

\section{Moduli of curves with a $G$-cover}\label{moduli1}

Consider $G$ a finite group, $\R_{g,G}$ is the moduli space of genus $g$ smooth curves
 with a principal $G$-bundle. 
 The moduli $\R_{g,G}$ comes with a natural forgetful proper morphism $\pi\colon\R_{g,G}\to\M_g$.
As shown in \cite{delmum69}, the moduli space $\overline \M_g$ of stable curves, is a compactification of $\M_g$.
In the case of principal $G$-bundles over stable curves, the
nodal singularities prevent the forgetful projection $\pi$
to be proper.

In order to find a compactification of $\R_{g,G}$ which is proper
over~$\overline\M_g$, we
introduce two equivalent stacks: the one of twisted $G$-covers of genus~$g$,
denoted by $\mc B^{\bal}_g(G)$, and the one of admissible $G$-covers of genus $g$, denoted by
$\adm_g^G$. These stacks are Deligne-Mumford and are proven to be isomorphic by Abramovich, Corti and Vistoli (see~\cite{acv03}),
 we introduce both of them because we will use different insights from both
points of view. 
The coarse space $\overline\R_{g,G}$ of these spaces is a compactification of
 $\R_{g,G}$, and it comes with a proper forgetful morphism
$\overline\R_{g,G}\to\overline\M_g$ which extends $\pi$.

\subsection{Curves with principal $G$-bundles} 
Given any finite group $G$, in this section we explore the
geometry of principal $G$-bundles over stable curves
and their automorphisms.

\subsubsection{Moduli of stable curves}
In \cite{delmum69}, Deligne and Mumford carry a local analysis of the
stack~$\overline\mm_{g,n}$ of stable curves, based on deformation theory. For every $n$-marked stable curve $(C;p_1,\dots,p_n)$,
the deformation functor is representable (see \cite{serne06} and \cite[\S11]{acgh11})
and it is represented by a smooth scheme $\Def(C;p_1,\dots,p_n)$ of
dimension $3g-3+n$ with one distinguished point $q$.
The deformation scheme comes with a universal family 
$X\to\Def(C;p_1,\dots,p_n)$
whose central fiber $X_{q}$ is identified with $(C;p_1,\dots,p_n)$. Every automorphism of the
central fiber naturally extends to the whole family $X$ by the universal
property of the deformation scheme.
The strict henselization
of $\overline\mm_{g,n}$ at the geometric point $[C;p_1,\dots,p_n]$ is the same of
the Deligne-Mumford stack
$$[\Def(C;p_1,\dots,p_n)\slash\Aut(C;p_1,\dots,p_n)]$$
at $q$. As a consequence,
for every geometric point $[C;p_1,\dots,p_n]$ of the coarse space~$\overline\M_{g,n}$, the strict
Henselization at $[C;p_1,\dots,p_n]$ is
$\Def(C;p_1,\dots,p_n)\slash\Aut(C;p_1,\dots,p_n)$.
This implies 
that every singularity of $\overline\M_{g,n}$ is a quotient singularity. From now on, we will refer
to the strict henselization of a scheme $X$ at a geometric point $q$
as the \cor{local picture} of $X$ at~$q$.

As showed in \cite[\S11.2]{acgh11}, given a smooth curve $C$
with $n$ marked points $p_1,\dots,p_n$, we have
$\Def(C;p_1,\dots,p_n)\cong H^1(C, T_C(-p_1-\dots-p_n))$,
where $T_C$ is the tangent bundle to curve $C$.

\begin{rmk}\label{rmk_sp}
Given a stable $n$-marked curve $C$, we denote by $C_1,C_2,\dots, C_V$
its irreducible components. Let $\nor\colon \overline C\to C$ be the normalization
morphism of $C$, and denote by $\overline C_i$ the normalization of component $C_i$ for
every $i$, then $\overline C=\sqcup_i\overline C_i$. We mark on $\overline C$
the preimage point via $\nor$ of any marked point or node.
We denote by $g_i$ the genus
of $\overline C_i$ for any $i$, by $D_i$ the divisor of marked points on $\overline C_i$ and by $n_i:=\deg (D_i)$ its degree.
The stability condition for $C$ is equivalent to  
$2g_i-2+n_i>0$ for all $i$.
\end{rmk}

\begin{rmk}\label{rmk_def1}
We follow \cite{chiofar12} to give a more explicit description
of the deformation scheme.
For a nodal curve $C$, consider
$\Def(C;\sing C)$ the universal deformation of curve $C$ alongside with
its nodes. We impose $n=0$ in this for sake of simplicity, the $n>0$ case
is similar. 
If $V$ is the number of irreducible components of $C$, there is a canonical decomposition
\begin{equation}
\Def(C;\sing C)=\bigoplus_{i=1}^V\Def(\overline C_i;D_i)\cong\bigoplus_{i=1}^VH^1(\overline C_i,T_{\overline C_i}(-D_i)).
\end{equation}
Furthermore, if $\delta$ is the number
of nodes of $C$, the quotient $\Def(C)\slash\Def(C;\sing C)$
has a canonical splitting 
\begin{equation}\label{eq_quot1}
\Def(C)\slash\Def(C;\sing C)=\bigoplus_{j=1}^\delta M_j,
\end{equation}
where $M_j\cong \A^1$ is the
 deformation scheme of node $q_j$ of $C$. The isomorphism $M_j\to \A^1$
is non-canonical and choosing one isomorphism is equivalent to choose a smoothing
of the node.\newline
\end{rmk}

\subsubsection{Group actions}\label{grthe}

Given any finite group $G$ and an element $h$ in it, we call
$c_h\colon G\to G$ the conjugation automorphism 
such that $c_h\colon g\mapsto h\cdot g\cdot h^{-1}$ for all $g$ in~$G$.
The subgroup of conjugation automorphisms, inside $\Aut(G)$, is called
group of the \cor{inner automorphisms} and  denoted by $\inn(G)$.
We call $\sub(G)$ the set of $G$ subgroups and, for any subgroup $H\in\sub (G)$,
we call $Z_G(H)$ its centralizer
$$Z_G(H):=\{g\in G|\ gh=hg\ \forall h\in H\}.$$
We denote by $Z_G$ the center of the whole group.
The group $\inn(G)$ acts naturally
on $\sub(G)$. 

\begin{defin}
We call $\mc T(G)$ the set of the orbits of the $\inn(G)$-action in $\sub(G)$. Equivalently,
$\mc T(G)$ is the set of conjugacy classes of $G$ subgroups.
\end{defin}
\begin{defin}\label{def_subc}
Consider two subgroup conjugacy classes $\mc H_1,\mc H_2$ in $ \mc T(G)$,
we say that $\mc H_2$ is a subclass of $\mc H_1$, denoted by $\mc H_2\leq \mc H_1$,
if for one element  $H_2\in \mc H_2$ (and hence for all), there exists $H_1\in\mc H_1$
such that $H_2$ is a subgroup of $H_1$. If the inclusion is strict, then $\mc H_2$ is a
strict subclass of $\mc H_1$ and the notation is $\mc H_2<\mc H_1$.
\end{defin}

Consider a transitive $G$-set $\mc T$, \cor{i.e.} a set $\mc T$ with a transitive left $G$-action 
$\psi\colon G\times \mc T\to \mc T$.
Any map $\eta\colon\mc T\to G$ induces, via $\psi$, a map $\mc T\to \mc T$. In particular,
$$E\mapsto \psi(\eta(E),E),\ \ \forall E\in \mc T.$$
This induces a map 
$$\psi_*\colon G^{\mc T}\to{\mc T}^{\mc T}.$$
If we denote by $S_{\mc T}\subset {\mc T}^{\mc T} $ the subset of $\mc T$ permutations,
we obtain that $\psi_*^{-1}(S_{\mc T})$ is the subset of
maps $\mc T\to G$ inducing a $\mc T$ permutation.

Consider an element $E$ in $\mc T$. We denote by $H_E$ its stabilizer,
\cor{i.e.}~the $G$ subgroup fixing~$E$. Given any other element
$\psi(g,E)$ for some $g$ in $G$, its stabilizer is 
$$H_{\psi(g, E)}=g\cdot H_E\cdot g^{-1},$$
this proves the following lemma.

\begin{lemma}\label{lemma_h}
Given any transitive $G$-set $\mc T$, there exists
a canonical conjugacy class $\mc H$ in~$\mc T(G)$, and a canonical surjection
$\mc T\twoheadrightarrow \mc H$
sending any element to its stabilizer.
\end{lemma}

Given the transitive $G$-set $\mc T$ and the group $G$ seen
as a $G$-set with respect to the $\inn(G)$-action,
we consider the set of
$G$-equivariant maps $\Ho^G(\mc T,G)$.

\begin{lemma}\label{lemma_h2}
For any element $E$ in $\mc T$, and any map $\eta$ in
$\Ho^G(\mc T,G)$, $\eta(E)\in Z_G(H_E)$.
\end{lemma}
\proof The equivariance condition means that
$$\eta(\psi(h,E))=c_h(\eta(E))=h\cdot \eta(E)\cdot h^{-1}$$
for all $h$ in $G$. If $h$ is in $H_E$, the left hand side
of the equality above is simply $\eta(E)$, therefore $c_h(\eta(E))=\eta(E)$
for all $h$ in $H_E$, and this is possible if and only if $\eta(E)$
is in $Z_G(H_E)$.\fine

\begin{prop}\label{prop_sdeng}
Given any object $E$ in $\mc T$, there exists a canonical isomorphism
$$\Ho^G(\mc T,G)\cong Z_G(H_E).$$
Moreover, the set of equivariant maps $\Ho^G(\mc T,G)$ is uniquely determined
by the canonical class $\mc H$ associated to $\mc T$ (see Lemma \ref{lemma_h}).
\end{prop}
\proof The first part of the proposition follows from Lemma \ref{lemma_h2}.
We observe that if we consider another object $E'=\psi(s,E)$,
then $H_{E'}=s\cdot H_E\cdot s^{-1}$ and 
$$Z_G(H_{E'})=s\cdot Z_G(H_E)\cdot s^{-1}.$$
Therefore there exists an inclusion $\Ho^G(\mc T,G)\hookrightarrow G$
which is determined, up to conjugation, by the class  $\mc H$ of $H_E$.\qed

\subsubsection{Principal $G$-bundles}

\begin{defin}[principal $G$-bundle]\label{pgb}
If $G$ is a finite group, a principal $G$-bundle over a scheme $X$ is
a fiber bundle $F\to X$ together with a left action $\psi\colon G\times F\to F$
such that the induced morphism
$$\widetilde \psi\colon G\times F\to F\times_X F,$$
is an isomorphism. Here 
 $\widetilde\psi=\psi\times \pi_2$, where $\pi_2$ is the projection $G\times F\to F$.
 \end{defin}
\begin{rmk}
As a direct consequence of the definition, every geometric fiber of $F\to X$ is isomorphic
to the group $G$ itself.
\end{rmk}
\begin{rmk}
The category of principal $G$-bundles is denoted by $BG$ and comes with
a natural forgetful functor $BG\to Sch$.
\end{rmk}

We introduce the stack $\rr_{g,G}$ of
smooth curves of genus $g$ with a principal $G$-bundle.

\begin{defin}\label{def_gbun}
In the category $\rr_{g,G}$, the objects are smooth $S$-curves $X\to S$ of genus~$g$,
equipped with a principal $G$-bundle $F\to X$, for any scheme $S$. 
The morphisms of $\rr_{g,G}$ are commutative diagrams 
$$
\xymatrix{
F'\ar[d]_{b}\ar[r] & X'\ar[d]\ar[r] &S'\ar[d]\\
F\ar[r] & X\ar[r] & S\\
}
$$
such that the two squares are cartesian and $b$ is $G$-equivariant
with respect to the natural $G$-actions.
The category $\rr_{g,G}$ comes with a forgetful functor
$\pi\colon \rr_{g,G}\to \mm_g$,
sending any object or morphism
on the underlying curve or curve morphism.
\end{defin}

Consider a connected normal scheme $X$
and a principal $G$-bundle $F\to X$. We denote by $\Aut_{\cov}(X,F)$
its  automorphism group in the category of coverings, that is the automorphisms 
of $F$ commuting with the projection $F\to X$. Furthermore, $\Aut_{BG}(X,F)$
is its automorphism group in the category of principal $G$-bundles, that
is the covering automorphisms of $F$ compatible with the
natural $G$-action.

We call $\mc T(F)$ the set of connected components of any
principal $G$-bundle $F\to X$. 
The group $G$ acts transitively on $\mc T(F)$,
and by abuse of notation we call 
$\psi\colon G\times \mc T(F)\to \mc T(F)$ this action.
As explained in Section \ref{grthe}, this action induces a map
$\psi_*\colon G^{\mc T(F)}\to {\mc T(F)}^{\mc T(F)}$.

\begin{prop}\label{prop_z}
If $X$ is a connected normal scheme, and $F\to X$ a
principal $G$-bundle, then we have the following canonical identifications:
\begin{itemize}
\item $\Aut_{\cov}(X,F)=\psi_*^{-1}(S_{\mc T(F)})$;
\item $\Aut_{BG}(X,F)=\Ho^G(\mc T(F),G)$.
\end{itemize}
Here we denoted by $S_{\mc T(F)}$ the set of $\mc T(F)$ permutations.
\end{prop}

\proof 
We start by showing the identification $\Ho_{\cov}(F\slash X,F\slash X)=G^{\mc T(F)}$,
where the first one is the set of covering automorphism of a principal $G$-bundle $F\to X$.
Consider any covering morphism 
$b\colon F\to F$ over $X$,
given the isomorphism 
$\widetilde\psi\colon G\times F\to F\times_X F$
introduced in Definition \ref{pgb}, we consider
the chain of maps
\begin{equation}\label{chain_eq}
F\xrightarrow{b\times \id} F\times_X F\xrightarrow{\widetilde\psi^{-1}}G\times F\xrightarrow{\pi_1} G.
\end{equation}
As $G$ is discrete, the map above is constant
on the connected components and therefore
it induces a map $\Ho_{\cov}(F\slash X,F\slash X)\xrightarrow{\eta} G^{\mc T(F)}$
which moreover is bijective.

The morphism $b$ is an automorphism if and only if $\eta(b)$
acts bijectively
on $\mc T(F)$, \cor{i.e.}~if and only if
$\psi_*(\eta(b))\in S_{\mc T(G)}$.\newline

The automorphisms of $F$ as a principal $G$-bundle must moreover
preserve the $G$-action, \cor{i.e.}~we must have 
$$b_h:=b\circ \psi(-,h)=h\cdot b(-)\ \ \forall h\in G,$$
where $\cdot$ is the multiplication in $G$. We observe that
$\eta(b_h)=(\eta(b)\circ\psi(-,h))\cdot h$, where
we denoted by $\psi$ the $G$-action on $F$ and $\mc T(F)$
indistinctly. Therefore $\eta(b)\circ\psi(-,h)=\eta(b_h)\cdot h^{-1}$, and so
$$\eta(b)\circ\psi(-,h)=c_h\circ\eta(b),$$
which is the exact definition of 
$\eta$ being in $\Ho^G(\mc T(F),G)\subset G^{\mc T(F)}$.\fine

\begin{rmk}\label{rmkconne}
In the case of a connected principal $G$-bundle $F\to X$, 
the proposition above summarizes in
$\Aut_{\cov}(X,F)=G$ and $\Aut_{BG}(X,F)=Z_G$.
\end{rmk}

For a general $G$-bundle $F\to X$,
the set of connected components $\mc T(F)$
has a transitive $G$-action. By Lemma \ref{lemma_h},
this induces a canonical conjugacy class $\mc H$ in $\mc T(G)$.

\begin{defin}
We call principal $\mc H$-bundle, a principal $G$-bundle whose
canonical associated class in $\mc T(G)$ is $\mc H$. Equivalently,
the stabilizer of every connected component in $\mc T(F)$ is
a $G$ subgroup in $\mc H$.
\end{defin}

\begin{rmk}
By Proposition \ref{prop_sdeng}, the automorphism group of
any principal $\mc H$-bundle, is isomorphic to $Z_G(H)$,
where $H$ is any $G$ subgroup in the $\mc H$ class.\newline
\end{rmk}

\subsection{Twisted $G$-covers}

To enlarge the notion of principal $G$-bundles we admit non-trivial stabilizers at the nodes
of a stable curve, by defining twisted curves. 
The twisting techniques are
widely discussed in \cite{abravis02} and 
\cite{acv03}, furthermore twisted curves
are introduced in \cite{chiofar12}
in the case of a level structure
on stable curves.

\subsubsection{Definitions}
\begin{defin}[Twisted curve]
A twisted $n$-marked $S$-curve is a diagram 
$$\begin{array}{ccc}
\Sigma_1,\Sigma_2,\dots,\Sigma_n & \subset &\ssC \\
 & & \downarrow \\
 & & C\\
 & & \downarrow \\
 & & S.
 \end{array}$$
 Where:
 \begin{enumerate}
 \item $\ssC$ is a Deligne-Mumford stack, proper over $S$,
 and étale locally it is a nodal curve over~$S$;
 \item the $\Sigma_i\subset \ssC$ are disjoint closed substacks in the smooth
 locus of $\ssC\to S$ for all~$i$;
 \item $\Sigma_i\to S$ is an étale gerbe for all~$i$;
 \item $\ssC\to C$ exhibits $C$ as the coarse space of $\ssC$, and it is an isomorphism over $C_{\gen}$.
 \end{enumerate}
\end{defin}

We recall that, given a scheme $U$ and a finite abelian group $\mmu$ acting on $U$,
the stack $[U\slash \mmu]$ is the category of principal $\mmu$-bundles $E\to T$, for
any scheme $T$, equipped with a $\mmu$-equivariant morphism $f\colon E\to U$.
The stack $[U\slash \mmu]$ is a proper Deligne-Mumford stack and has a natural
morphism to its coarse scheme $U\slash\mmu$.\newline

By the definition of twisted curve we get the local pictures:
\begin{itemize}
\item \cor{at a marking}, morphism $\ssC\to  C\to S$ is locally isomorphic to 
$$\left[\Spec A[x']\slash \mmu_r\right]\to \Spec A[x]\to\Spec A$$
for some normal ring $A$ and some integer $r>0$.
Here $x=(x')^{r}$, and
$\mmu_r$ is the cyclic group of order $r$ acting on $\Spec A[x']$ by the action 
$\xi\colon x'\mapsto \xi x'$ for any $\xi\in\mmu_r$;\newline

\item \cor{at a node}, morphism $\ssC \to C\to  S$ is locally isomorphic to 
$$\left[\Spec\left( \frac{A[x',y']}{(x'y'-a)}\right)\slash \mmu_r\right]\to \Spec\left( \frac{A[x,y]}{(xy-a^\ell)}\right)\to \Spec A$$
for some  integer $r>0$ and $a\in A$. Here $x=(x')^\ell,\ y=(y')^\ell$. The group $\mmu_r$ acts
by the action 
$$\xi\colon(x',y')\mapsto (\xi x',\xi^m y')$$ where $m$ is an
element of $\Z\slash r$ and $\xi$ is a primitive $r$th root of the unit. The action is called \cor{balanced}
if $m\equiv -1 \mod r$. A curve with balanced action at every node is called a balanced curve.
\end{itemize}

\begin{defin}[Twisted $G$-cover]
Given an $n$-marked twisted curve 
$(\Sigma_1,\dots,\Sigma_n;\ \ssC\to C\to S)$,
 a twisted $G$-cover is a representable 
stack morphism $\phi\colon\ssC\to BG$,
\cor{i.e.}~an object of the category $\Fun(\ssC,BG)$
which moreover is representable.
\end{defin}
\begin{defin}
We introduce category 
 $\mc B_{g,n}(G)$. 
Objects of $\mc B_{g,n}(G)$ are twisted $n$-marked $S$-curves of genus $g$ with a twisted
$G$-cover, for any scheme $S$. 

Consider two twisted $G$-covers $\phi'\colon \ssC'\to BG$
and $\phi\colon \ssC\to BG$ over the twisted $n$-marked curves $\ssC'$ and $\ssC$ respectively.
A morphism $(\ssC',\phi')\to(\ssC,\phi)$
is a pair $(\msf f, \alpha)$
such that $\msf f\colon \ssC'\to \ssC$
is a morphism of $n$-marked twisted curves,
and $\alpha\colon \phi'\to \phi\circ\msf f$
is an isomorphism in $\Fun(\ssC',BG)$.
\end{defin}

Following \cite{abravis02}, the category $\mc B_{g,n}(G)$
can be defined as the $2$-category of twisted stable
$n$-pointed maps of genus $g$ and degree $0$ to the category $BG$.
In the same paper it is observed that the automorphism group of every $1$-morphism
is trivial, therefore this $2$-category is equivalent to the category obtained by replacing
$1$-morphisms with their $2$-isomorphism classes. In \cite{abravis02} this category is denoted by
$\mc K_{g,n}(BG,0)$, the notation $\mc B_{g,n}(G)$ for the case of twisted $G$-covers
appears for example in \cite{acv03}.

\begin{defin}
A balanced twisted $G$-cover is a twisted $G$-cover
over a twisted balanced curve. We call $\mc B_{g,n}^{\bal}(G)$ the full
sub-functor of twisted balanced $G$-covers.
\end{defin}

Twisted $G$-covers generalize the notion of root of the trivial bundle.
Indeed, for any twisted curve $\ssC$ and any integer $\ell>0$, there exists a canonical bijection
between the set of twisted $\mmu_\ell$-covers over $\ssC$, and
the set of $\ell$th roots of $\Oo_{\ssC}$.
Here a faithful line bundle is a 
line bundle $\ssL\to \ssC$ such that the associated morphism $\ssC\to B\C^*$
is representable,
and an $\ell$th root of $\Oo_{\ssC}$
is a faithful line bundle such that $\ssL^{\xx\ell}\cong\Oo_{\ssC}$.\newline

\subsubsection{Local structure of twisted covers}\label{sec_locsttwcov}
We consider a twisted curve $\ssC$ over a geometric point $\Spec(\C)$.
For any marked or nodal point $p$, 
the local picture of $\ssC$ at $p$ is the same as $[U\slash\mmu_r]$ at the origin,
for some scheme $U$ and positive integer $r$.
Any principal $G$-bundle over~$\ssC$, or equivalently
any object of $BG(\ssC)$, is locally isomorphic at $p$ to a principal $G$-bundle on~$[U\slash\mmu_r]$.

In \cite[\S2.1.8]{acv03} is explained how to realize twisted stable maps
as twisted objects over scheme theoretic curves. In particular, a principal $G$-bundle on $[U\slash\mmu_r]$ 
is the same as a principal $G$-bundle $f\colon \tilde F\to U$ with the natural $G$-action $\psi\colon G\times \tilde F\to \tilde F$,
and also with a $\mmu_r$-action
$\nu\colon \mmu_r\times \tilde F\to \tilde F$ which is compatible
with the $\mmu_r$-action on $U$ and with $\psi$.

In formulas we have:
\begin{enumerate}
\item $f\circ \nu(\xi,-)=\xi\cdot f\colon \tilde F\to U$, for all $\xi\in\mmu_r$;
\item $\psi(h,\nu(\xi,-))=\nu(\xi,\psi(h,-))\colon \tilde F\to \tilde F$, for all $h\in G$ and~$\xi\in\mmu_r$.
\end{enumerate}

\begin{rmk}\label{rmk_0}
We consider at first the case of a marked point $p$ of $\ssC$ whose
local picture is $[\A^1\slash\mmu_r]$ with $\mmu_r$ acting by multiplication. By what we just said
we have a principal $G$-bundle $\tilde F\to \A^1$, and for any 
$\xi\in\mmu_r$ a morphism $\tilde\alpha(\xi)\colon \tilde F\to\tilde F$
such that
$$\tilde\alpha(\xi):=\nu(\xi,-).$$
If we fix a privileged $r$th root $\xi_r=\exp(2\pi\slash r)$, then
$\tilde \alpha(\xi_r)(\tilde p)=\psi(h_{\tilde p},\tilde p)$,
for all preimages $\tilde p$ of~$p$, where $ h_{\tilde p}$ is an element of group $G$ depending on $\tilde p$.
\end{rmk}

\begin{rmk}\label{rmk_1}
In the case of a node $q$ of $\ssC$, its local picture is $[V\slash \mmu_r]$
for some positive integer $r$ where $V\cong\{x'y'=0\}\subset\A^2_{x',y'}$ and the $\mmu_r$-action
is given by $\xi\cdot (x',y')=(\xi x',\xi^{-1}y')$ for all $\xi\in \mmu_r$.

The normalization of the node neighborhood $V$ is naturally isomorphic to
$\A^1_{x'}\sqcup \A^1_{y'}\to V$.
We consider the normalization $\nor\colon\overline\ssC\to \ssC$ of the twisted curve $\ssC$,
the local picture of $\nor$ morphism at $q$ is
$$[\A_{x'}^1\slash\mmu_{r_q}]\sqcup [\A_{y'}^1\slash\mmu_{r_q}]\to [V\slash\mmu_{r_q}].$$
We denote by $q_1\in \A^1_{x'}$ and $q_2\in \A^1_{y'}$ the two preimages of $q$.
As before, a twisted $G$-cover on $[V\slash\mmu_r]$ is the same as a principal $G$-bundle $\tilde F\to V$
plus a $\mmu_r$-action on $V$ with the right compatibilities. This induces
\begin{itemize}
\item two principal $G$-bundles 
$\tilde F'\to \A^1_{x'}$ and $\tilde F''\to \A^1_{y'}$
with the naturally associated $\mmu_{r_q}$-actions.
We denote by $\nu'\colon\mmu_{r_q}\times\tilde F'\to \tilde F'$ and
$\nu''\colon\mmu_{r_q}\times \tilde F''\to \tilde F''$ these actions;
\item a gluing isomorphism
between the central fibers
$\kappa_q\colon \tilde F'_q\xrightarrow{\cong} \tilde F''_q$.
This means that 
\begin{enumerate}[label=\roman{*}., ref=(\roman{*})]
\item $\kappa_q(\psi(h,-))=\psi(h,\kappa_q(-))\colon \tilde F'_q\to \tilde F''_q$ for any $h\in G$;
\item $\kappa_q(\nu'(\xi,-))=\nu''(\xi,\kappa_q(-))\colon \tilde F'_q\to \tilde F''_q$ for any $\xi\in\mmu_{r_q}$.
\end{enumerate}
And furthermore, $\tilde F=(\tilde F'\sqcup \tilde F'')\slash \kappa_q$.
\end{itemize}

Following Remark \ref{rmk_0}, we define
$\alpha'(\xi):=\nu'(\xi,-)\colon F'\to F'$,
$\alpha''(\xi):=\nu''(\xi,-)\colon F''\to F''$
for any $\xi\in\mmu_{r_q}$.
By the balancing condition, if we have two points
$\tilde q_1$ and $\tilde q_2$ in $F'_{q_1}$ and $F''_{q_2}$ respectively,
such that $\kappa_q(\tilde q_1)=\tilde q_2$, then
$h_{\tilde q_1}=h_{\tilde q_2}^{-1}$.\newline
\end{rmk}

This local structure can be encoded in conjugation classes associated to every marked or nodal point.
Consider a marked point $p$ of $\ssC$ and the local picture $[\A^1\slash \mmu_r]$
at $p$, then the twisted $G$-cover $\phi\colon\ssC\to BG$ induces a morphism
$\phi\colon[\A^1\slash\mmu_r]\to BG$.
This induces a morphism $\tilde\phi_p\colon \mmu_r\to G$
defined up to conjugation,
which  is an injection by the $\phi$ representability.

\begin{defin}\label{defin_gt0}
The conjugacy class $\lq \tilde \phi_p\rq $ of $\tilde\phi_p$ is called $G$-type of $\phi$ at $p$.
\end{defin}

In the case of a node $q$, the composition of $\phi$ with the normalization induces
$$\tilde\phi_{q_1}\colon\mmu_r\to G\ \ \m{and}\ \ \tilde\phi_{q_2}\colon\mmu_r\to G,$$
and by the balancing condition the two $G$-types are the inverse of each other,
$\lq \tilde\phi_{q_1}\rq=\lq \tilde\phi_{q_2}\rq^{-1}$.
Once we choose a privileged branch of a node, we call $G$-type
of that node the $G$-type with respect to the
restriction of the cover to that branch. Switching the branches changes the $G$-type
into its inverse class.\newline

\subsubsection{Local structure of $\mc B^{\bal}_{g,n}(G)$}

The local structure of $\mc B^{\bal}_{g,n}(G)$ can be described
with a very similar approach to what we did for $\overline\mm_{g,n}$. We work
the case $n=0$ of unmarked twisted $G$-covers. Given
a twisted $G$-cover $(\ssC,\phi)$, its deformation functor is representable and
the associated scheme $\Def(\ssC,\phi)$ is isomorphic to $\Def(\ssC)$
via the  forgetful functor $(\ssC,\phi)\mapsto \ssC$.
The automorphism group $\Aut(\ssC,\phi)$ naturally
acts on $\defo(\ssC,\phi)=\defo(\ssC)$ and 
the local picture of $\mc B^{\bal}_g(G)$ at $[\ssC,\phi]$
is the same of  $[\defo(\ssC)\slash\Aut(\ssC,\phi)]$
at the central point.

\begin{rmk}\label{rmk_def2}
Consider a twisted curve $\ssC$ whose coarse
space is the curve $C$, we give a description of the scheme $\defo(\ssC)$ as
we did in Remark \ref{rmk_def1} for $\defo(C)$ with the notation of Remark \ref{rmk_sp}.
As $\ssC$ is a twisted curve, every node $q_j$ has a possibly non-trivial
stabilizer, which is a cyclic group of order $r_j$.

The deformation $\defo(\ssC;\sing\ssC)$ of $\ssC$ alongside with its nodes,
is canonically identified with the deformation of $C$ alongside with its nodes
$\defo(\ssC;\sing \ssC)=\defo(C;\sing C)$.
As in the previous case, the following quotient has a canonical splitting.
\begin{equation}\label{eq_quot2}
\defo(\ssC)\slash \defo(\ssC;\sing \ssC)=\bigoplus_{j=1}^{\delta}R_j.
\end{equation}
In this case $R_j\cong \A^1$ is the deformation scheme of the
node $q_j$ together with its stack structure. If we consider the schemes $M_j$
of Equation (\ref{eq_quot1}) in Remark \ref{rmk_def1}, there exists
for every $j$ a canonical morphism $R_j\to M_j$ of order $r_j$ ramified in exactly one point.\newline
\end{rmk}

\subsection{Admissible $G$-covers}
In order to define admissible $G$-covers, in the next sections we introduce admissible covers and
we put a balancing condition on them.

\subsubsection{Admissible covers}

\begin{defin}[Admissible cover]
Given a nodal $S$-curve $X\to S$ with marked points,
an admissible cover $u\colon F\to X$ is a morphism such that:
\begin{enumerate}
\item the composition $F\to S$ is a nodal $S$-curve;
\item given a geometric point $\bar s \in S$, every node of $F_{\bar s}$ maps via $u$ to a node of $X_{\bar s}$;
\item the restriction $F|_{X_{\gen}}\to X_{\gen}$ is an étale cover of degree $d$;
\item given a geometric point $\bar s \in S$, the local picture of $F_{\bar s}\xrightarrow{u} X_{\bar s}$ at a point of $F_{\bar s}$ mapping to a marked point of $X$
is isomorphic to
$$\Spec A[x']\to \Spec A[x]\to \Spec A,$$
for some normal ring $A$, an integer $r>0$ and $u^*x=(x')^r$;
\item the local picture of $F_{\bar s}\xrightarrow{u} X_{\bar s}$ at a node of $F_{\bar s}$
is isomorphic to
$$\Spec \left(\frac{A[x',y']}{(x'y'-a)}\right)\to \Spec \left(\frac{A[x,y]}{(xy-a^r)}\right)\to \Spec A,$$
for some integer $r>0$ and an element $a\in A$, $u^*x=(x')^r$ and $u^*y=(y')^r$.
\end{enumerate}
\end{defin}

The category $\adm_{g,n,d}$ of $n$-pointed stable curves of genus $g$
with an admissible cover of degree $d$, is a proper Deligne-Mumford stack.\newline

Consider $F\to C$ an admissible cover of a nodal curve $C$, a $G$-action on
$F$ such that the restriction $F|_{C_{\gen}}\to C_{\gen}$ is a principal $G$-bundle,
 a smooth point $p$ of $C$, and a preimage $\tilde p\in F$ of $p$. We denote by $H_{\tilde p}\subset G$
the stabilizer of $\tilde p$. By definition of admissible cover, if $p\in C_{\gen}$ (\cor{i.e.}~$p$ is non-marked), then
$H_{\tilde p}=(1)$. Moreover, the $G$-action induces a primitive character
$$\chi_{\tilde p}\colon H_{\tilde p}\to \GL(T_{\tilde p}F)=\C^*,$$
where $T_{\tilde p}F$ is the tangent space to $F$ in $\tilde p$.

Given any subgroup $H\subset G$, for any primitive character $\chi\colon H\to \C^*$ and for any $s\in G$,
we denote by $\chi^s$ the conjugated character $\chi^s\colon sHs^{-1}\to \C^*$
such that $\chi^s(h)=\chi(s^{-1}hs)$ for all $h\in G$.

In the set of pairs $(H,\chi)$, with $H$ a $G$ subgroup and
$\chi\colon H\to \C^*$ a character, we introduce
the equivalence relation 
$(H,\chi)\sim (H',\chi')$ iff there exists $s\in G$ such that $H'=sHs^{-1}$ and $\chi'=\chi^s$.
Consider a point $\tilde p$ on $F$ with stabilizer $H_{\tilde p}$
and associated character $\chi_{\tilde p}$. We observe that for any point $s\cdot \tilde p$ of the same fiber,
$$H_{s\cdot \tilde p}=s H_{\tilde p}s^{-1}\ \ \m{and}\ \ \chi_{s\cdot \tilde p}=\chi_{\tilde p}^s.$$
Therefore the equivalence class of the pair $(H_{\tilde p},\chi_{\tilde p})$
only depends on the point $p$.

\begin{defin}\label{def_locgtype}
For any smooth point $\tilde p$ on $F$, we call 
\cor{local index} the associated pair $(H_{\tilde p},\chi_{\tilde p})$.
For any smooth point $p\in C$,
the conjugacy class of the local index of any $\tilde p$ in $F_p$ 
is called the \cor{$G$-type} at~$p$, following the notation in \cite{berroma11}.
We denote the $G$-type by
$\lq H_p,\chi_p\rq $,
where $H_p$ is the stabilizer of one of the points in $F_p$, and $\chi_p$ the associated character.
\end{defin}

The notion of $G$-type is equivalent to the one introduced in Definition \ref{defin_gt0}.
We will discuss this equivalence in \S\ref{sec_equiv}.

\subsubsection{Balancing the $G$-action}
\begin{lemma}\label{lem_stab}
Consider $u\colon F\to C$ an admissible cover of a nodal curve $C$
such that the restriction $F|_{C_{\gen}}\to C_{\gen}$ is a principal $G$-bundle.
If $\tilde p\in F$ is one of the preimages of a node or a marked point, then the stabilizer $H_{\tilde p}$
is a cyclic group.
\end{lemma}
\proof If $\tilde p$ is the preimage of a marked point, the local picture of morphism $u$ 
at $\tilde p$
is 
$$\Spec A[x']\to \Spec A[x],$$
where $x'=x^r$ for some integer $r>0$. This local description
induces an action of $H_{\tilde p}$ on 
$U:=\Spec A[x']$ which is free
and transitive on $U\backslash \{\tilde p\}$. The group
of automorphisms of $U\backslash\{\tilde p\}$ preserving $r$ is
exactly $\mmu_	r$, therefore $H_{\tilde p}$ must be cyclic too.

In the case of a node $\tilde p$ we observe that $u$ is locally
isomorphic to
$$\Spec \left(\frac{A[x',y']}{(x'y'-a^r)}\right)\to\Spec \left(\frac{A[x,y]}{(xy-a)}\right),$$
where $x'=x^r$ and $y'=y^r$, for an integer $r>0$ and an element $a\in A$.
The scheme $U':=\Spec \left(A[x',y']\slash(x'y'-a^r)\right)$ is the union
of two irreducible components $U_1,U_2$, and 
we can apply the deduction above to $U_i\backslash\{\tilde p\}$ for $i=1,2$.\fine

Observe that the set of characters $\chi\colon \mmu_r\to \C^*$ of a cyclic group, 
is the group $\Z\slash r\Z$.
In particular, the character associated to $k\in \Z\slash r\Z$ maps $\xi\mapsto \xi^k$ for any $\xi$
$r$th root of the unit.

Focusing on the case of a node $\tilde p\in F$, we observe that $H_{\tilde p}$ acts independently on
the two branches $U_1$ and $U_2$. We denote by $\chi_{\tilde p}^{(1)}$ and $\chi_{\tilde p}^{(2)}$
the characters of these actions.

\begin{defin}
The $G$-action at node $\tilde p$ is balanced when $\chi_{\tilde p}^{(1)}=-\chi_{\tilde p}^{(2)}$,
that is they are opposite as elements of $\Z\slash r\Z$ (where $r$ depends on the $\tilde p$ fiber).
\end{defin}

\begin{defin}[Admissible $G$-cover]\label{def_admcov}
Given any finite group $G$, consider an admissible cover $F\to C$ of a nodal curve $C$, it is an admissible $G$-cover if
\begin{enumerate}
\item the restriction $F|_{C_{\gen}}\to C_{\gen}$ is a principal $G$-bundle.
This implies,
by Lemma \ref{lem_stab}, that for every node or marked point $\tilde p\in F$, the stabilizer $H_{\tilde p}$ is a cyclic group;
\item the action of $H_{\tilde p}$ is balanced for every node $\tilde p\in F$.
\end{enumerate}
\end{defin}

This notion was firstly developed
by Abramovich, Corti and Vistoli in~\cite{acv03}, and also
by Jarvis, Kaufmann and Kimura in~\cite{jkk05}.

\begin{defin}
We call $\adm_{g,n}^G$ the stack of stable curves of genus $g$ with $n$ marked points
and equipped with an admissible $G$-cover.
\end{defin}

\begin{rmk}\label{rmk_locind}
For any cyclic subgroup $H\subset G$, the image of a 
 character $\chi\colon H\to \C^*$
is the group of $|H|$th roots of the unit, $\mmu_{|H|}$. We choose a privileged
root in this set, which is $\exp(2\pi i\slash |H|)$. After this choice, 
The datum of $(H,\chi)$, is equivalent to the datum of the $H$ generator
$h=\chi^{-1}(e^{2\pi i\slash |H|})$. As a consequence, the conjugacy class $\lq H,\chi\rq $
is identified with the conjugacy class $\lq h\rq $ of $h$ in~$G$.
\end{rmk}

\begin{defin}\label{def_hd}
Given an admissible $G$-cover $F\to C$ over an $n$-marked stable curve, the series 
$\lq h_1\rq ,\ \lq h_2\rq ,\dots,\ \lq h_n\rq$,
of the $G$-types of the singular fibers over the marked points,
is called \cor{Hurwitz datum} of the cover.
The stack of admissible $G$-covers of genus $g$ with a given
Hurwitz datum is denoted by $\adm_{g,\lq h_1\rq ,\dots,\lq h_n\rq }^G$.
\end{defin}

\begin{rmk}\label{def_gtn}
Given an admissible $G$-cover $F\to C$, if $p$ is a node of $C$
and $\tilde p$ one of its preimages on $F$, then 
the local index of $\tilde p$ and the $G$-type of $p$
are well defined once we fix a privileged branch of $p$.
Switching the branches sends the local index and the $G$-type
in their inverses.
\end{rmk}

Consider a smooth curve $C$ of genus $g$ and $n$ marked points $p_1,\dots, p_n$,
the fundamental group of $C_{\gen}=C\backslash\{p_1,\dots,p_n\}$
has $2g+n$ generators
$\alpha_1,\alpha_2,\dots\alpha_g,\beta_1\dots,\beta_{g},\gamma_1,\dots,\gamma_n$.
These generators respect the following relation,
\begin{equation}\label{eq_boh}
\alpha_1\beta_1\alpha_1^{-1}\beta_1^{-1}\cdots\alpha_{g}\beta_{g}\alpha_{g}^{-1}\beta_{g}^{-1}\cdot\gamma_1\cdots\gamma_n=1,
\end{equation}
and this is sufficient to represent the fundamental group. This is called the canonical representation
of the fundamental group of a genus $g$ smooth curve.

It is possible to describe admissible $G$-covers over smooth curves
by the monodromy action, as done for example in \cite[\S 2.3]{berroma11} and \cite[\S 3.5]{schmizel18}.
Consider a smooth curve $C$, a generic point $p_*$ on it
and an admissible $G$-cover $F\to C$. We denote the points of the fiber $F_{p_*}$ by
$\tilde p_*^{(g)}$ for any $g\in G$, in such a way that $g\cdot \tilde p_*^{(1)}=\tilde p_*^{(g)}$.
This induces a group morphism $\pi_1(C_{\gen},p_*)\to G$.
This monodromy morphism is well defined
up to relabelling  the points $\tilde p_*^{(i)}$, \cor{i.e.}~up to
$G$ conjugation.
The following proposition is a rephrasing of \cite[Lemma 2.6]{berroma11}.

\begin{prop}\label{prop_gfg2}
Given a smooth $n$-marked curve $(C;p_1,\dots,p_n)$
and a point $p_*$ on its generic locus $C_{\gen}$,
the set of isomorphism classes of admissible $G$-covers on $C$
 is naturally in bijection
with the set of conjugacy classes of maps 
$$\varpi\colon \pi_1(C_{\gen},p_*)\to G.$$
\end{prop}
\begin{rmk}\label{rmk_mono}
We also point out that the monodromy of $\gamma_i$
 at any point~$p_*^{(g)}$, with $g\in G$, is given by a small circular
lacet around the deleted point~$p_i$. Therefore
by definition of $G$-type, if $\lq h_i\rq $ is the $G$-type of $p_i$,
then $\lq \varpi(\gamma_i)\rq=\lq h_i\rq$.\newline
\end{rmk}

\subsubsection{Admissible $G$-cover automorphisms}\label{admaut}

Consider an admissible $G$-cover $F\to C$ over a smooth
$n$-marked curve $(C;p_1,\dots,p_n)$.
We denote by $\mc T(F)$ the set of connected
components of $F$, which inherits the $G$-action $\psi$.
For any connected component $E\subset F$, we denote by $H_E\subset G$
its stabilizer.
The component $\psi(s, E)$, for some element $s$ of $G$, has
stabilizer $s\cdot H_E\cdot s^{-1}$.
Therefore the conjugacy class of the stabilizer is independent 
on the choice of $E$.  As in the case of principal $G$-bundles,
for every admissible $G$-cover there exists a canonical class
$\mc H$ in $\mc T(G)$ such that the stabilizer of every $E$ in $\mc T(F)$
is a subgroup $H_E$ in $\mc H$. Moreover, we have a canonical surjective map
$$\mc T(F)\twoheadrightarrow \mc H.$$
\begin{defin}
Given the set $\mc T(G)$ of subgroup conjugacy classes in $G$,
and a class $\mc H$ in it,
an admissible $\mc H$-cover
is an admissible $G$-cover such that every connected component has stabilizer in $\mc H$.
\end{defin}

\begin{defin}
We denote by $\adm^{G,\mc H}_g$ the 
stack of admissible $\mc H$-covers over stable curves
of genus $g$, and we denote by
$\adm^{G,\mc H}_{g,\lq h_1\rq,\dots,\lq h_n\rq}$ the stack
of admissible $\mc H$-cover with Hurwitz datum $\lq h_1\rq,\dots,\lq h_n\rq$
over the $n$ marked points.
\end{defin}

It is possible to generalize
the second point of Proposition \ref{prop_z}. We denote by $\Aut_{\adm}(C,F)$ the
set of automorphisms of an admissible $G$-cover $F\to C$. 

\begin{prop}\label{prop_zz}
Consider $(C;p_1,\dots,p_n)$ a nodal $n$-marked curve, and
$F\to C$ an admissible $G$-cover, then
$$\Aut_{\adm}(C,F)=\Ho^G(\mc T(F), G).$$
\end{prop} 
\proof In the case of a smooth curve $C$, we consider
the general locus $C_{\gen}=C\backslash\{p_1,\dots,p_n\}$.
The restriction $F|_{C_{\gen}}$ is a principal $G$-bundle, 
therefore by Proposition \ref{prop_z},
$$\Aut_{\adm}(C,F)\subset \Aut_{BG}(C_{\gen},F_{\gen})=\Ho^G(\mc T(F|_{C_{\gen}}),C_{\gen}).$$
Since $\mc T(F|_{C_{\gen}})=\mc T(F)$ and every automorphism 
of $F|_{C_{\gen}}\to C_{\gen}$ extends to the whole $F$,
the thesis follows in this case.

In the case of a general stable curve $C$, with
$C_1,\dots, C_V$ its connected components, and 
$F_i$ the restriction $F|_{C_i}$ for any $i$, as a consequence of the first part, we have
$$\Aut_{\adm}(C_i,F_i)=\Ho^G(\mc T(F_i),C_i).$$ 
The balancing condition at the nodes imposes that any automorphism
in $\Aut_{\adm}(C,F)$ acts as the same multiplicative factor
on two touching components. This means that 
a sequence of functions in 
$\prod_i\Ho^G(\mc T(F_i),G)$, induces a global automorphism
if and only if it is the sequence of restrictions of
a global function $\Ho^G(\mc T(F),G)$.\fine

\subsubsection{Equivalence between twisted and admissible $G$-covers}\label{sec_equiv}
We introduced the two categories $\mc B_g^{\bal}(G)$ and $\adm_g^G$ with
the purpose of ``well'' defining the notion of principal $G$-bundle over stable
non-smooth curves. These two categories are proven isomorphic in \cite{acv03}.

\begin{teo}[see {\cite[Theorem 4.3.2]{acv03}}]\label{teo_3e}
There exists a base preserving equivalence between $\mc B_g^{\bal}(G)$ and $\adm_g^G$,
therefore in particular they are isomorphic Deligne-Mumford stacks.
\end{teo}

The proof proposed in \cite{acv03} can be sketched quickly. Given a twisted $G$-cover 
$\phi\colon \ssC\to BG$, the restriction to $\ssC_{\gen}=C_{\gen}$
is a principal $G$-bundle $F_{\gen}\to C_{\gen}$ on the generic locus of the coarse space $C$,
and this can be uniquely completed to an admissible $G$-cover $F\to C$.
Conversely, given an admissible $G$-cover $F\to C$, it induces a quotient stack
$\ssC:=[F\slash G]$ and therefore a representable morphism $\ssC\to BG$
with balanced action on nodes.

In what follows we will adopt the notation $\overline\rr_{g,G}$ for the equivalent stacks
$\mc B_g^{\bal}(G)$ and $\adm_g^G$.
For every class $\mc H$ in $\mc T(G)$ we denote by $\overline\rr_{g,G}^{\mc H}$
the full substack of $\overline\rr_{g,G}$ whose objects are admissible $\mc H$-covers.

The correspondence of Theorem \ref{teo_3e} allows the translation of every machinery we developed on 
twisted $G$-covers to admissible $G$-covers, and conversely.
For example, the two definitions of $G$-type we introduced are equivalent.
Precisely, consider a twisted $G$-cover $(\ssC,\phi)$, a point $p$ whose
$G$-type is $\lq\tilde\phi_p\rq$, and an element $\tilde\phi_p\colon\mmu_r\to G$ in the
class of the $G$-type. Therefore $\im\tilde\phi_p=H\subset G$ is a cyclic subgroup and
$\tilde\phi_p^{-1}\colon H\to \mmu_r$ is a character. The class $\lq H,\tilde\phi_p^{-1}\rq$
is precisely the $G$-type at $p$ from the admissible $G$-cover point of view.

Furthermore, we can use over twisted $G$-covers the notion of Hurwitz datum.  We 
will denote by 
$\overline\rr_{g,\lq h_1\rq ,\dots,\lq h_n\rq }^{\mc H}$
the stack of admissible $\mc H$-covers of genus $g$ with Hurwitz datum $\lq h_1\rq ,\dots,\lq h_n\rq $.

If there is no risk of confusion, we will say that a twisted $G$-cover $(\ssC,\phi)$
``is'' an admissible $G$-cover $F\to C$ (or the other way around), meaning that $F\to C$ is
the naturally associated admissible $G$-cover to $(\ssC,\phi)$.\newline

\section{Dual graphs and torsors}\label{chap_mod}
In this section we introduce the important tool of dual graphs to 
describe subloci of the moduli of curves with a twisted $G$-cover. This subject
was already treated by the author in \cite{gale15} in the case
of spin curves. Here we update this tool in order to generalize this notion
to the case of $G$-covers. Furthermore, we introduce  torsors 
and some of their fundamental properties.

\subsection{Decorated dual graphs and $G$-covers}

\subsubsection{Basic graph theory}\label{sec_graph}
Consider a graph $\Gamma$ with vertex set $V$ and edge set~$E$, we call
\cor{loop} an edge that starts and ends on the same vertex, we call
\cor{separating} an edge $e$ such that the graph with vertex set $V$ and
edge set $E\backslash \{ e\}$ is disconnected. 
We denote by $E_{\sep}$ the set
of separating edges, and by $\E$ the set of oriented
edges: the elements of this set are edges in $E$ equipped with an orientation, 
in particular for every edge $e\in\E$ we denote by
$e_+$ the head vertex and by $e_-$ the tail. There is a canonical $2$-to-$1$ projection $\E\to E$.
We also introduce a conjugation 
in $\E$, such that for each $e\in \E$,  the conjugated edge $\bar e$ is obtained
by reversing the orientation, in particular
$(\bar e)_+=e_-$.
For every graph $\Gamma$, 
when there is no risk of confusion we denote by $V$
the cardinality of the vertex set $V(\Gamma)$ and by $E$
the cardinality of the edge set~$E(\Gamma)$.\newline

We consider a finite group $G$ acting on graph $\Gamma$:
That is, we consider two $G$-actions
on the vertex set and on the edge set,
$$G\times V(\Gamma)\to V(\Gamma)\ \ \ \m{and}\ \ \ G\times \E(\Gamma)\to \E(\Gamma).$$
We denote these actions by $h\cdot v$ and $h\cdot e$ for every
$h$ in $G$ and every vertex $v$ and oriented edge $e$.
These actions must 
respect the following natural intersection conditions
\begin{enumerate}
\item $(h\cdot e)_+=h\cdot e_+\ \ \forall h\in G, e\in \E(\Gamma)$;
\item $\overline{h\cdot e}=h\cdot \overline e\ \ \forall h\in G, e\in \E(\Gamma)$.
\end{enumerate}
Observe that there are no faithfulness conditions, therefore
any vertex or edge may have a non-trivial
stabilizer. We denote by $H_v$ and $H_e$ the stabilizers of vertex
$v$ and edge $e$ respectively. We have
$H_{s\cdot v}=s\cdot H_v\cdot s^{-1}$ for any $v\in V(\Gamma)$ and $s\in G$,
and the same is true for $H_e$. In general, every orbit of vertices (or oriented edges)
is characterized by a conjugacy class $\mc H$ in~$\mc T(G)$, and every
element of $\mc H$ is the stabilizer of some object in the orbit.

\begin{defin}[Cochains]
The group of $0$-cochains is the group of $G$-valued functions on $V(\Gamma)$
compatible with the $G$-action
$$C^0(\Gamma;G):=\left\{a\colon V(\Gamma)\to G|\ a(g\cdot v)=g\cdot a(v)\cdot g^{-1}\right\}.$$

The group of $1$-cochains is the group of antisymmetric functions on $\E$
with the same compatibility condition
$$C^1(\Gamma;G):=\left\{b\colon \E\to G|\ b(\bar e)=b(e)^{-1},\ b(g\cdot e)=g\cdot b(e)\cdot g^{-1}\right\}.$$
\end{defin}

These groups generalize the cochain groups defined by
Chiodo and Farkas in \cite{chiofar12}. In particular the Chiodo-Farkas groups
refer to the case of a trivial $G$-action on $\Gamma$.\newline

There exists a natural differential $\delta\colon C^0(\Gamma;G)\to C^1(\Gamma;G)$
such that
$$\delta a(e):=a(e_+)\cdot a(e_-)^{-1},\ \ \forall a\in C^0(\Gamma;G)\ \forall e\in \E.$$
Consider the set $\mc T(\Gamma)$ of the connected components of the graph,
with the naturally induced $G$-action.
The exterior differential fits into an useful exact sequence of groups
\begin{equation}\label{seq_1}
0\to \Ho^G(\mc T(\Gamma),G)\xrightarrow{i} C^0(\Gamma;G)\xrightarrow{\delta} C^1(\Gamma;G).
\end{equation}
Here the injection $i$ sends $f\in \Ho^G(\mc T(\Gamma),G)$ on the cochain $a$
such that for every component $\gamma\in\mc T(\Gamma)$, 
$a$ is constantly equal to $f(\gamma)$ on $\gamma$. If $\Gamma$ is a connected graph,
then the first term of the exact sequence is the group $G$ and $i$ sends $g\in G$
to the associated constant cochain.
We recall that for any group we can define a (non-associative) $\Z$-action via 
$h\cdot n:=h^n$ for all $h$ in $G$ and $n\in \Z$.

\begin{prop}\label{prop_imdelta.1}
A $1$-cochain $b$ is in $\im\delta$ 
if and only if, for every circuit $\mc K=(e_1,\dots, e_k)$ in $\E$, we have
$$b(\mc K) :=b(e_1)\cdot b(e_2)\cdots b(e_k)=1.$$
\end{prop}
\proof If $b\in \im\delta$, the condition above is easily verified. To complete
the proof we show that if the condition si verified,
then there exists a cochain $a\in C^0(\Gamma;G)$ such that
$\delta a=b$. We choose a vertex $v\in V(\Gamma)$ and impose
$a(v)=1$, for any other vertex $w\in V(\Gamma)$ we consider a path 
$\mc P=(e_1,\dots,e_m)$ starting in $v$ and ending in $w$.
We set 
$$a(w):=b(\mc P)=b(e_1)\cdots b(e_m).$$
By the condition on circuits, the cochain
$a$ is well defined, independently of path $\mc P$, and by construction we have $b=\delta a$.\fine

\subsubsection{Trees and tree-like graphs}\label{treetree}
\begin{defin}
A tree is a graph that does not contain any circuit.
A \cor{tree-like} graph is a connected graph whose only circuits are loops.
\end{defin}
\begin{rmk}\label{ineq_graph}
For every connected graph $\Gamma$, 
the first Betti number  $b_1(\Gamma)=E-V+1$
is the dimension
rank of the homology group $H_1(\Gamma;\Z)$. 
Note that, $b_1$ being positive,
$E\geq V-1$.
This inequality is an equality if and only if $\Gamma$ is a tree.
\end{rmk}

For every connected graph $\Gamma$ with vertex set $V$ and
edge set $E$, we can choose a connected subgraph $T$ with the same vertex set
and edge set $E_T\subset E$ such that $T$ is a tree.
\begin{defin}
 The graph
$T$ is called a \cor{spanning tree} for $\Gamma$.
\end{defin}

\begin{lemma}\label{ineq_graph2}
If $E_{\sep}\subset E$ is the set of edges in $\Gamma$ that are
separating, then 
$E_{\sep}\leq V-1$
 with equality if and only if $\Gamma$
is  tree-like.
\end{lemma}
\proof
If $T$ is a spanning tree for $\Gamma$ and $E_T$ its edge set, then
$E_{\sep}\subset E_T$. Indeed, an edge $e\in E_{\sep}$ is the only path between its two extremities, therefore,
since $T$ is connected, $e$ must be in $E_T$.
Thus $E_{\sep}\leq E_T=V-1$, with equality if and only if all the edges
of $\Gamma$ are loops or separating edges, \cor{i.e.}~if $\Gamma$ is a tree-like graph.\fine

\subsubsection{Graph contraction and graph $G$-covers}
Consider a graph $\Gamma$ with vertex set $V$ and
edge set $E$, we choose a subset $D\subset E$
which is stable
by the $G$-action.
\begin{defin}
Consider the graph $\Gamma_0$ such that:
\begin{enumerate}
\item the edge set of $\Gamma_0$ is $E_0:=E\backslash D$;
\item given the relation in $V$, $v\sim w$ if $v$ and $w$ are linked by an edge $e\in D$, the
vertex set of $\Gamma_0$ is $V_0:=V\slash\sim$.
\end{enumerate}
The graph $\Gamma_0$ inherits naturally a $G$-action. The natural morphism $\Gamma\to\Gamma_0$
is called \cor{contraction} of $D$ or $D$-contraction.
\end{defin}

Edge contraction will be useful, in particular
we will consider the image of the exterior differential $\delta$
and its restriction over contractions of a given graph.
If $\Gamma_0$ is a contraction of $\Gamma$, 
then $E(\Gamma_0)$ is canonically a subset of $E(\Gamma)$.
As a consequence, cochains over $\Gamma_0$
are cochains over $\Gamma$ with the additional condition that
the values on $E(\Gamma)\backslash E(\Gamma_0)$ are all the identity. Then
we have a natural immersion
$C^i(\Gamma_0;G)\hookrightarrow C^i(\Gamma;G)$.
 Consider the two exterior differentials
$$\delta\colon C^0(\Gamma;G)\to C^1(\Gamma;G)\ \ \ \m{and}\ \ \ \delta_0\colon C^0(\Gamma_0;G)\to C^1(\Gamma_0;G).$$
The following proposition follows.
\begin{prop}\label{prop_imcap}
The differential $\delta_0$ is the restriction of $\delta$ on $C^0(\Gamma_0;G)$.
$$\im\delta_0=C^1(\Gamma_0;\Z\slash\ell)\cap \im\delta.$$
\end{prop}

Given any graph $\Gamma$ with a $G$-action, we define its $G$-quotient
$\Gamma \slash G$
by $V(\Gamma\slash G):=V( \Gamma)\slash G$ and $E(\Gamma\slash G):=E(\Gamma)\slash G$. The conditions
on the $G$-action assure
that $\Gamma\slash G$ is well defined. Moreover, the edge contraction of a subset $D\subset E( \Gamma)$ stable under $G$-action,
is compatible with the quotient, so that if $ \Gamma\to \Gamma_0$ is the $D$-contraction,
then 
$\Gamma\slash G\to \Gamma_0\slash G$
is the contraction of $D\slash G$ (the $G$-action on the new quotiented graphs is trivial).

We call a $G$-graph morphism $\tilde\Gamma\to \Gamma$ a graph $G$-cover if $\Gamma\cong \tilde\Gamma\slash G$
and $\tilde\Gamma\to\Gamma$ 
is the natural quotient morphism.
For any vertex $ \tilde v$ of $ \tilde\Gamma$, we denote by $H_{\tilde v}$ its stabilizer
in $G$.
For any vertex $v$ of~$\Gamma$, its preimages in $V(\tilde\Gamma)$ all have
a stabilizer in the same conjugacy class $\mc H$ in~$\mc T(G)$, \cor{i.e.}~for 
all $\tilde v$ in $f^{-1}(v)$ we have $H_{\tilde v}\in \mc H$. 
Moreover, for every
subgroup $H$ in the class $\mc H$, there exists a vertex preimage $\tilde v$  of $v$
with stabilizer exactly $ H$. In particular the cardinality of the $v$ fiber
is $|G|\slash |H| $ where $|H|$ is the cardinality of any subgroup in $\mc H$.
The same is true for any edge $e$ in $E(\Gamma)$.

We observe that it is possible to give another description of the cochain
groups of $\tilde \Gamma$ by considering the graph $G$-cover $\tilde \Gamma\to \Gamma$.
Given a set $T$ with a $G$-action and a conjugation $e\mapsto \bar e$, define $\overline\Ho^G(T,G)$
as the set of morphisms $f\colon T\to G$ compatible with the $G$-action and such that $f(\bar e)=f(e)^{-1}$.
We extend the $G$-action on $\E(\tilde \Gamma)$ by defining it as a fibered product in the category
of $G$-sets,
$\E(\tilde \Gamma):=E(\tilde\Gamma)\times_{E(\Gamma)}\E(\Gamma)$,
this prevents that $h\cdot e=\bar e$ for some $e\in \E(\tilde\Gamma)$ and $h\in G$.

\begin{prop}\label{prop_coch}
Consider a graph $G$-cover $f\colon\tilde\Gamma\to \Gamma$. We have the identification
$$C^0(\tilde\Gamma;G)=\prod_{v\in V(\Gamma)}\Ho^G(f^{-1}(v),G).$$
Moreover, 
$$C^1(\tilde\Gamma;G)=\prod_{e\in \E(\Gamma)}\overline\Ho^G(f^{-1}(e),G).$$\newline
\end{prop}

\subsubsection{Graph $G$-cover of an admissible $G$-cover}

\begin{defin}[dual graph]
Consider a nodal curve $C$, its dual graph $\Gamma(C)$
is defined~by
\begin{align*}
V(\Gamma(C)) &:=\{\m{irreducible components of }C\}\\
E(\Gamma(C)) &:=\{\m{nodes of }C\}
\end{align*}
with the natural link relations. 
\end{defin}
\begin{rmk}
We observe that the set
of oriented edges $\E(\Gamma(C))$ is naturally identified with the
set of nodes equipped with a privileged branch, or equivalently with the set of
node preimages on the normalization $\overline C$.
\end{rmk}

For any admissible $G$-cover $F\to C$, consider
the dual graphs $\tilde \Gamma:=\Gamma(F)$ and $\Gamma:=\Gamma(C)$.
Therefore $\Gamma=\tilde\Gamma\slash G$
and $\tilde\Gamma\to \Gamma$ is a graph $G$-cover.
We recall the correspondence between admissible $G$-covers over a stable curve $C$,
and twisted $G$-covers over $C$, treated in section~\ref{sec_equiv}. As a consequence,
the dual graphs $\tilde \Gamma$ and $\Gamma$ introduced for any
admissible $G$-cover, are well defined for the
associated twisted $G$-cover, too.

Consider the function $b_F$ defined on $\E(\tilde \Gamma)$ that sends
 any oriented edge $\tilde e$ to the local index $(H,\chi)$
 of the associated node, where the privileged node branch (necessary
 to define the local index) is given by the $\tilde e$ orientation (see the
 Definition \ref{def_gtn} of the local index).
We observe that for any $h\in G$, $b_F(h\cdot \tilde e)=(h Hh^{-1},\chi^h)$.
Furthermore, we associate to $b_F$ another function $M_{b_F}$
sending any $e$ in $\E(\Gamma)$ to the $G$-type $\lq H,\chi\rq$
of the associated node.

\begin{defin}\label{def_tfac}
We call \cor{index cochain} of the admissible $G$-cover $F\to C$, the
function~$b_F$. Moreover, we call
type function of $F\to C$, the function
$M_{b_F}$.
When there is no risk of confusion, we denote the type function
of $F\to C$ simply by $M$.
\end{defin}

\begin{rmk}
Once we choose a privileged $r$th root of the unit $\xi_r=\exp(2\pi i\slash r)$
for every positive integer $r$, the index cochain $b_F$ is identified
with a $1$-cochain in $C^1(\tilde \Gamma; G)$,
and the associated type function is a function
$M_{b_F}\colon \E(\Gamma)\to \lq G\rq $. 
\end{rmk}

\begin{rmk}
In the case of an admissible $G$-cover
with $G$ abelian group, the type function uniquely
determines the index cochain. 
In the case of $G=\mmu_\ell$, our notation reduces
to the 
 multiplicity index notation of Chiodo and Farkas \cite{chiofar12}.
\end{rmk}

We observe that
the order of $M_{b_F}(e)$ is well defined for any $e\in \E(\Gamma)$ as the order of any element in the conjugacy class,
therefore we define the function $r\colon \E(\Gamma(C))\to \Z_{>0}$. Clearly $r(e)=r(\bar e)$ for any $e$.

\begin{defin}\label{contr_graph}
A pair $(\Gamma,r(-))$, where $r\colon \E(\Gamma)\to \Z_{>0}$
is an even function, is called
\emph{decorated graph}.
The pair $(\Gamma(C),r(-))$ given by the admissible $G$-cover $F\to C$ (or
equivalently, the associated twisted $G$-cover $(\ssC,\phi)$) is called decorated graph of the cover.
If there is no risk of confusion, we will refer also to $\Gamma(C)$ or $\Gamma$ alone
as the decorated graph.
\end{defin}

Let $D\subset \E(\tilde \Gamma)$ be the subset of edges where the cochain $b_F$ of local indices is trivial, that~is
$$D:=\{\tilde e\in \E(\tilde\Gamma)|\ b_F(\tilde e)=1\}.$$

\begin{defin}\label{gamma0}
The graph $\tilde \Gamma_0$ is the result of the $D$-contraction on $\tilde \Gamma$. The graph $\Gamma_0$ is the quotient
$\tilde\Gamma_0\slash G$. Equivalently, it is the graph $\Gamma$ after the conctraction
of the edges where the type function $M$ has value $\lq 1\rq $.\newline
\end{defin}

\subsection{Basic theory of sheaves in groups and torsors}\label{sec_shea}
In this section we refer in particular
to Calmès and Fasel paper \cite{calfas15} for notations and definitions.
Consider a scheme $S$ and a site $\mb T$ over the category $Sch\slash S$ of $S$-schemes.
An $S$-sheaf for us will be a sheaf over $(Sch\slash S,\mb T)$.
Consider $\mc G$ an $S$-sheaf in groups, and $P$ an $S$-sheaf
in sets with a left $\mc G$-action.

\begin{defin}[torsor]
The sheaf $P$ is a torsor under $\mc G$, or a $\mc G$-torsor, if 
\begin{enumerate}
\item the application 
$\mc G\times P\to P\times P$,
where the components are the action and the identity, is an isomorphism;
\item for every covering $\{S_i\}$ of $S$, $P(S_i)$ is non-empty for every $i$.
\end{enumerate}
\end{defin}

For example, if $G$ is a finite group, a principal $G$-bundle over a scheme $S$, is 
a $\mc G$-torsor, where $\mc G$ is the $S$-sheaf in groups
defined by $\mc G(S'):=S'\times G$ for any $S$-scheme $S'$.
When we consider any $S$-sheaf in groups $\mc G$ as acting
on itself, we get a $\mc G$-torsor called
\cor{trivial $\mc G$-torsor}.

\begin{prop}[see {\cite[Proposition 2.2.2.4]{calfas15}}]
An $S$-sheaf $P$ with a left $\mc G$-action is a torsor
if and only if it is $\mb T$-locally isomorphic to the trivial
torsor $\mc G$.
\end{prop}

Consider two $S$-sheaves $P$ and $P'$ with $\mc G$-action respectively
on the left and on the right. 

\begin{defin}\label{def_cp}
We denote by $P'\wedge^{\mc G}P$ the cokernel sheaf
of the two morphisms 
$$\mc G\times P'\times P\rightrightarrows P'\times P$$
given by the $\mc G$-action on $P$ and $P'$ respectively.
This is called \cor{contracted product}.
Equivalently, $P'\wedge^{\mc G}P$ is the sheafification of
the presheaf of the orbits of $\mc G$ acting on $P'\times P$
by
$$(h,(z',z))\mapsto (z'h^{-1},hz).$$
\end{defin}
\begin{rmk}
If $\mc G$ is the sheaf in groups constantly equal to $\C^*$
and $P,P'$ are two line bundles, then the contracted product
is simply the usual tensor product $P\xx P'$.
\end{rmk}

If another $S$-sheaf in groups $\mc G'$ acts
on the left on $P'$, then the contracted product
$P'\we^{\mc G} P$ has a $\mc G'$-action on the left, too.
The same is true for a $\mc G'$-action on the right on $P$.

\begin{lemma}[see {\cite[Lemma 2.2.2.10]{calfas15}}]
The $\wedge$ construction is associative. Consider $\mc G$
and $\mc G'$ two $S$-sheaves in groups,
$P$ and $P'$ two $S$-sheaves with respectively left $\mc G$-action and
right $\mc G'$-action, finally $P''$ an $S$-sheaf with $\mc G'$-action
on the left and $\mc G$-action on the right, and the actions commute.
Then there exists a canonical isomorphism
$$(P'\wedge^{\mc G'}P'')\wedge^{\mc G} P\cong P'\we^{\mc G'}(P''\we^{\mc G}P).$$
Moreover, we have $\mc G\we^{\mc G} P\cong P$ for every $\mc G$-torsor $P$.
\end{lemma}

\begin{prop}[see {\cite[Proposition 2.2.2.12]{calfas15}}]\label{prop_pprod}
Consider a morphism $\mc G\to \mc G'$ of $S$-sheaves
in groups and the associated $\mc G$-action (on the right) on $\mc G'$.
The map
$$P\mapsto \mc G'\we^{\mc G} P,$$
from the category of $\mc G$-torsors to $\mc G'$-torsors,
is a functor.
\end{prop}

\begin{defin}
Given an $S$-scheme $S'$ and a site $\mb T$ on $Sch\slash S$, we denote by $H^1_{\mb T}(S',\mc G)$
the pointed set of $\mc G$-torsors (on the left) over $S'$ with respect
to the $\mb T$ topology. The base point of the set being
the torsor $\mc G$ itself.
\end{defin}

We observe that if $P'$ is a $\mc G$-bitorsor, on the left and on the right, over $S'$, then the 
contracted product $P'\we^{\mc G}P$ is a $\mc G$-torsor (on the left)
for every $\mc G$-torsor $P$. Therefore $P'$ induces a map
$$P'\we^{\mc G}-\colon H^1_{\mb T}(S',\mc G)\to H^1_{\mb T}(S',\mc G).$$

This cohomology type notation fits with the cohomology type behavior
we are going to describe.
We refer for the following results 
 to \cite[\S 2.2.5]{calfas15} or \cite[Chap.3]{gir71}. Consider three $S$-sheaves
in groups fitting in a short exact sequence
\begin{equation}\label{eq_sestor}
1\to \mc G_1\to \mc G_2\to \mc G_3\to 1.
\end{equation}

\begin{teo}\label{teo_torsors}
This gives a long exact sequence in cohomology 
\begin{equation}\label{eq_lestor}
1\to \mc G_1(S)\to \mc G_2(S)\xrightarrow{\delta} \mc G_3(S)\xrightarrow{\tau} H^1_{\mb T}(S,\mc G_1)\xrightarrow{w} H^1_{\mb T}(S,\mc G_2)\to H^1_{\mb T}(S,\mc G_3).
\end{equation}
This is an exact sequence of pointed sets, and it is exact in $\mc G_1(S)$ and $\mc G_2(S)$ as a sequence
of groups.
\end{teo}

To describe the map $\tau$, observe that $\mc G_3=\mc G_2\slash\mc G_1$.
By \cite[Proposition 3.1.2]{gir71}, the set $\mc G_3(S)$ is in bijection
with the set of sub-$\mc G_1$-torsors of $\mc G_2$. Any
object $Q$ in $\mc G_3(S)$ in sent by $\tau$ on the $\mc G_1$-torsor
induced by the pullback along $\mc G_2\to\mc G_2\slash\mc G_1$.
As a consequence $\tau(Q)$ is a $\mc G_1$-bitorsor.  

Via the $\tau$ map we also have a $\mc G_3(S)$-action on $H^1_{\mb T}(S,\mc G_1)$. Indeed,
for every $Q$ in $\mc G_3(S)$ and for every $\mc G_1$-torsor $P$, we obtain by
contracted product the $\mc G_1$-torsor
$\tau(Q)\we^{\mc G_1} P$.\newline

To state the next proposition, we observe that $\mc G_1$ acts trivially on the right on $\mc G_3$,
therefore given any $\mc G_1$-torsor $P$ (on the left), we have
the identification of sheaves 
$\mc G_3\we^{\mc G_1}P=\mc G_3$.
Consider the map $\mc G_2\to\mc G_3$ in the short exact sequence (\ref{eq_sestor}), and its image via the
contracted product functor of Proposition \ref{prop_pprod},
$$\mc G_2\we^{\mc G_1} P\xrightarrow{-\we P}\mc G_3\we^{\mc G_1}P=\mc G_3.$$
We define $\mc G_2^P:=\mc G_2\we^{\mc G_1} P$, and 
$\delta^P\colon \mc G_2^P(S)\to \mc G_3(S)$.

\begin{prop}[see {\cite[Proposition 3.3.3]{gir71}}]\label{prop_fond}
For every $P$ in $H^1_{\mb T}(S,\mc G_1)$,
the stabilizer of $P$ with respect to the $\mc G_3(S)$-action induced by $\tau$, is the image
of $\delta^P\colon\mc G_2^P(S)\to \mc G_3(S)$.\newline
\end{prop}

\section{Singularities of $\overline{\mathcal R}_{g,G}$}\label{soR}
\subsection{Ghost automorphisms of a twisted curve}\label{sec_ghosta}

Consider a twisted $G$-cover $(\ssC,\phi)$, 
its automorphism group is
$$\Aut(\ssC,\phi):=\{(\msf f,\rho)|\ \msf f\in \Aut(\ssC),\ \rho\colon\phi\xrightarrow{\cong}\msf f^*\phi\}.$$

We observe that this group does not act faithfully on the universal
deformation $\defo(\ssC,\phi)$. Indeed, Proposition \ref{prop_zz}
describes the group $\Aut_{\ssC}(\ssC,\phi)$ of automorphisms of $(\ssC,\phi)$ acting trivially on $\ssC$, and these
automorphisms are the ones acting trivially on $\defo(\ssC,\phi)$, too. 
It becomes natural to
consider the group
$$\underline\Aut(\ssC,\phi):=\Aut(\ssC,\phi)\slash\Aut_{\ssC}(\ssC,\phi)=\{\msf f\in \Aut(\ssC)|\ \msf f^*\phi\cong\phi \m{ as twisted}\ G\m{-covers}\}.$$
\begin{rmk}\label{rmk_local}
The local description of $\overline \R_{g,G}$ at $[\ssC,\phi]$ could be rewritten $$\defo(\ssC)\slash\underline\Aut(\ssC,\phi).$$
\end{rmk}

The coarsening $\ssC\to C$ induces moreover a group morphism
$\underline\Aut(\ssC,\phi)\to \Aut(C)$.
We denote the kernel and the image of this morphism by
 $\underline\Aut_C(\ssC,\phi)$ and $\Aut'(C)$ (see also \cite[chap. 2]{chiofar12}).
They fit into the following short exact sequence,
\begin{equation}\label{seq_aut}
1\to \underline\Aut_C(\ssC,\phi)\to\underline\Aut(\ssC,\phi)\to \Aut'(C)\to 1.
\end{equation}
\begin{defin}
The group $\underline\Aut_C(\ssC,\phi)$ is called the group of ghost
automorphisms of $(\ssC,\phi)$.
\end{defin}

To describe the ghost automorphisms of a twisted $G$-cover, we start by describing $\Aut_C(\ssC)$,
the group of ghost automorphisms of the curve (not necessarily lifting to the cover).
Consider a node $q$ of $\ssC$ whose local picture is $[\{x'y'=0\}\slash\mmu_r]$.
Given an automorphism $\eta\in\Aut_C(\ssC)$,
the local action of $\eta$ at $q$ can be represented by an automorphism of $V=\{x'y'=0\}\subset\A^2$
such that
$$(x',y')\mapsto(\xi x',y'),$$
with $\xi$ a primitive root in $\mmu_r$.
We observe moreover that $(\xi x',y')\equiv (\xi^{u+1}x',\xi^{-u}y')$ for any integer $u$,
by the $\mmu_r$-action on $V$. Anyway, when it is not specified otherwise, 
we will use the lifting that acts trivially on the $y'$ coordinate.
Consider the dual decorated graph $(\Gamma(C),r(-))$
associated to the twisted $G$-cover $(\ssC,\phi)$, by definition
$r(e)$ is the order of the $q$-stabilizer where $q$ is the node associated
to edge $e\in E(\Gamma(C))$. We naturally extend the 
function $r$  over $\E(\Gamma)$.
As a consequence of the definition of $\Aut_C(\ssC)$, the action
of $\eta$ outside the $\ssC$ nodes is trivial. 
Then the whole
group $\Aut_C(\ssC)$ is generated by automorphisms of the form
$(x',y')\mapsto (\xi x', y')$ on a node, and trivial elsewhere.
We are interested in representing $\Aut_C(\ssC)$ as acting on the edges of
the dual graph, thus we introduce the following group.
\begin{defin}
Consider a decorated graph $(\Gamma,r(-))$, 
we denote by $r_{\lcm}$
 the least common multiple of all the orders $r(e)$ of the edges of $\Gamma$.
We define the group
$$S(\Gamma;r(-)):=\{f\colon \E(\Gamma)\to \Z\slash r_{\lcm}|\ f(e)=f(\bar e)\in \Z\slash r(e)\subset \Z\slash r_{\lcm}\}.$$
\end{defin}
We recall that $E(\Gamma)$ is the set of $\Gamma$ edges while 
$\E(\Gamma)$ is the set of $\Gamma$ edges with an orientation.
If $e\in \E(\Gamma)$ is an oriented $\Gamma$ edge, then we denote by
$\bar e$ the same edge with reversed orientation.

If $(\Gamma(\ssC),r(-))$ is the decorated dual graph
of the twisted $G$-cover $(\ssC,\phi)$, we define a morphism
$ S(\Gamma(\ssC);r(-))\to \Aut_C(\ssC),$ sending any function $\msf a$
on the automorphism whose action at the node associated to $e\in E(\Gamma)$ is 
$$(x',y')\mapsto (\msf a(e)\cdot x',y').$$
The morphism above is a canonical isomorphism, and we have the following identification
$$\underline\Aut_C(\ssC)= S(\Gamma(\ssC);r(-))= \bigoplus_{e\in E(\Gamma)}\mmu_{r(e)}.$$
Clearly the action is trivial on nodes with order $r=1$, so $S(\Gamma(\ssC);r(-))=S(\Gamma_0(\ssC);r(-))$.
We observe again that by choosing a privileged $r$th root $\exp(2\pi i\slash r)$ for any positive
integer~$r$, $\mmu_r$ is identified to $\Z\slash r$, and then $S(\Gamma_0(\ssC);r(-))\cong \bigoplus_{e\in E(\Gamma_0)} \Z\slash r(e)$.\newline

The group of ghost automorphisms
$$\underline\Aut_C(\ssC,\phi)=\{\msf a\in \Aut_C(\ssC)|\ \msf a^*\phi\cong\phi\}$$
 is a subset of $\Aut_C(\ssC)$. To describe it
 we will characterize the automorphisms in $\Aut_C(\ssC)$ lifting to the twisted $G$-cover~$\phi$.
 
If $C$ is the coarse space of $\ssC$, we consider the admissible $G$-cover $F\to C$
associated to $(\ssC,\phi)$,
and the normalization morphism $\nor\colon \overline C\to C$. We denote by $C_i$ the irreducible
components of $C$, by $\overline C_i$ their normalizations and by $F_i:=F|_{C_i}$ the $F$ restrictions. For any open subscheme $U\hookrightarrow C$,
$F|_U\to U$ is an admissible $G$-cover. Finally we define the pullbacks
$\overline F:=\nor^* F$ and $\overline U:= \nor^* U$. 
Consider the category $Sch\slash C$ of $C$-schemes and the Zarisky site $\mb T_{\Zar}$ on it.
Given the definition of automorphisms for
admissibile $G$-covers as stated in \S\ref{admaut}, we introduce the following definition.

\begin{defin}
The $C$-sheaf in groups $\H^F$ is defined for any open $C$-scheme $U\hookrightarrow C$ by,
$$\H^F(U):=\Aut_{\adm}(U,F|_U).$$
\end{defin}

We observe that $F$ is a $C$-sheaf
with a left $\H^F$-action, and we have a short exact sequence of $C$-sheaves in groups,
\begin{equation}\label{ses_nor}
1\to \H^F\to \nor_*\nor^*\H^F\xrightarrow{t}\H^F|_{\sing C}\to 1.
\end{equation}
The central sheaf is defined over any open subscheme $U\hookrightarrow C$ as 
$$\nor_*\nor^*\H^F(U)=\Aut_{\adm}(\overline U,\overline F|_{\overline U}).$$
There exists a $2:1$ cover $\overline F|_{\sing C}\to F|_{\sing C}$.
If $\varepsilon$ is a section of $\nor_*\nor^*\H^F(U)$,
its image via $t$ is obtained on every point $p$ of $F|_{\sing}(U)$
by taking the difference between the actions of $\varepsilon$ on the two preimages,
and therefore $t(\varepsilon)$ is well defined up to
ordering the branches of every node.

We pass to the associated long exact sequence. 
We observe that $\H^F(C)=\Ho^G(\mc T(F),G)=\Ho^G(\mc T(\tilde\Gamma),G)$ by Proposition \ref{prop_z}.
Moreover,
$$\nor_*\nor^*\H^F(C)=\Aut_{\adm}(\overline C,\overline F)=\prod_i\Ho^G(\mc T(F_i),G),$$
because the $\overline C_i$ are the connected components of $\overline C$.
If we denote by $f\colon\tilde \Gamma\to \Gamma$ the graph $G$-cover
associated to $F\to C$, by Proposition \ref{prop_coch} we have
$$\nor_*\nor^*\H^F(C)=C^0(\tilde \Gamma;G).$$
Finally, if $q_1,\dots,q_{\delta}$ are the nodes of $C$, then
$$\H^F|_{\sing C}(C)=\prod_j\Aut_{\adm}(q_j,F_{q_j})=\prod_j\Ho^G(F_{q_j},G).$$
By the definition of the dual graph $\tilde \Gamma$, the right hand side
of the equality above is identified with $\prod_{e\in E(\Gamma)}\overline\Ho^G(f^{-1}(e),G)$,
and by Proposition \ref{prop_coch} we have
$\H^F|_{\sing C}(C)\cong C^1(\tilde\Gamma;G)$.\newline

We consider the long exact sequence (\ref{eq_lestor}) associated to any short exact
sequence of $C$-sheaves. In our case the site over $Sch\slash C$ is the Zariski site $\mb T_{\Zar}$
and taking the long exact sequence associated to (\ref{ses_nor}), we get 
\begin{equation}\label{les_nor}
1\to \Ho^G(\mc T(\tilde \Gamma),G)\xrightarrow{i} C^0(\tilde\Gamma;G)\xrightarrow{\delta} C^1(\tilde\Gamma;G)\xrightarrow{\tau} H^1_{\mb T_{\Zar}}(C;\H^F)\xrightarrow{w}
H^1_{\mb T_{\Zar}}(C;\nor_*\nor^*\H^F)\to 1,
\end{equation}
where $H^1_{\mb T_{\Zar}}(C;\H^F)$ is the set of $\H^F$-torsors,
$H^1_{\mb T_{\Zar}}(C;\nor_*\nor^*\H^F)$ is the set of $\nor_*\nor^*\H^F$-torsors on $C$, 
and it is identified with $H^1_{\mb T_{\Zar}}(\overline C;\nor^*\H^F)$.
Moreover, the only object of $H^1_{\mb T_{\Zar}}(C;\H^F|_{\sing})$ is the trivial torsor.
The first part of this sequence is exactly the sequence~(\ref{seq_1}).
To describe explicitly the map $w$, consider
the normalization $\nor\colon \overline C\to C$. Then,
$$w\colon (F\to C)\mapsto (\overline F=\nor^*F\to \overline C).$$

Given any cochain $b\in C^1(\tilde\Gamma;G)$, by what we saw in Section \ref{sec_shea}
we know that $\H^F$ acts on the right on $\tau(b)$.
Therefore we can define an admissible $G$-cover by  the contracted product 
$\tau(b)\we^{\H^F} F$ (see Definition \ref{def_cp}).

Recall that to every admissible $G$-cover $F\to C$ is assigned an index cochain $b_F$ (see Definition \ref{def_tfac}).
Now consider an automorphism $\msf a\in \Aut_C(\ssC)=S(\Gamma(\ssC);r(-))$.
We define a $1$-cochain $b_F\cdot \msf a\in C^1(\tilde\Gamma(\ssC); G)$:
for every oriented edge $\tilde e$ of $\tilde \Gamma(\ssC)$, $b_F(\tilde e)=(H,\chi)$
that is a character $\chi\colon H\to \C$ where $H\subset G$ is a cyclic subgroup.
If $e\in \E(\Gamma(\ssC))$ is the projection of $\tilde e$, we define
 $$(b_F\cdot \msf a) (\tilde e):=\chi^{-1}(\msf a(e))\in H\subset G.$$

\begin{prop}\label{prop_pb}
Given a finite group $G$ and a twisted $G$-cover $(\ssC,\phi)$
consider the associated admissible $G$-cover $F\to C$,
whose index cochain is $b_F$. If $\msf a\in \underline\Aut_C(\ssC)=S(\Gamma(\ssC),r(-))$
is a ghost automorphism of $\ssC$, the pullback twisted $G$-cover $(\ssC,\msf a^*\phi)$,
where $\msf a^*\phi=\phi\circ\msf a$, is associated to the admissible $G$-cover
$$ \tau(b_F\cdot\msf a)\we^{\H^F}F.$$
\end{prop}

\proof Consider a node $q$ of the twisted curve $\ssC$. In Remark \ref{rmk_1}
we observed that the local picture of $(\ssC,\phi)$ at $q$ can be seen as a twisted object on $V\cong\{x'y'=0\}$.
This is equivalent to a principal $G$-bundle $\tilde F\to V$ with a compatible $\mmu_r$-action and the other
conditions of the same Remark. We remark that $\tilde F\slash \mmu_r\to V\slash\mmu_r\cong V$ is isomorphic to the local picture
of $F\to C$ around $q$. 

We start by characterizing $\msf a^*\phi$ with
respect to $\phi$. Again from Remark \ref{rmk_1} we know that 
$\tilde F=(\tilde F'\sqcup \tilde F'')\slash \kappa_q$
locally at node $q$, where $\tilde F'$ and $\tilde F''$ are the
two pullbacks of $\tilde F$ on the node branches $\A^1_{x'},\ \A^1_{y'}$,
and $\kappa_q\colon \tilde F'_q\to \tilde F''_q$
is the gluing morphism of the central fibers.
We consider the oriented edge $e\in \E(\Gamma(\ssC))$ associated
to $q$ with privileged branch $\A^1_{x'}$. We can lift the $\msf a$-action to the twisted object by
acting trivially on $\A^1_{y'}$ and by $\msf a(e)$ multiplication on~$\A^1_{x'}$.
By the same Remark \ref{rmk_1} we can lift the action to $\tilde F'$,
$$
\begin{tikzcd}
\tilde F'\ar[rr, "\alpha'(\msf a(e))"]\ar[d] && \tilde F'\ar[d]\\
\A^1_{x'}\ar[rr, "\msf a(e)\cdot -"] &&\A^1_{x'}.
\end{tikzcd}
$$
We observe that $\msf a^*\tilde F'\cong \tilde F'$ and $\msf a^*\tilde F''\cong \tilde F''$,
what really changes is the gluing morphism. Indeed,
$$\msf a^*\tilde F=(\tilde F'\sqcup\tilde F'')\slash(\kappa_q\circ\alpha'(\msf a(e))).$$
By definition of $\alpha'$, for any point $\tilde q'$ on the fiber $\tilde F'_q$,
we have $\alpha'(\msf a(e))(\tilde q')=\nu'(\msf a(e),\tilde q')$. Again by the $\alpha'$ definition,
$\alpha'(\msf a(e))(\tilde q')=\psi(h_{\tilde q'},\tilde q')$, where $(H,\chi)$
is the local index at $\tilde q'$ and $h_{\tilde q'}=\chi^{-1}(\msf a(e))$, that is $h_{\tilde q'}=(b_F\cdot \msf a)(\tilde e)$,
where $\tilde e\in \E(\tilde \Gamma)$ is the edge associated to $\tilde q'$ and
the privileged branch associated to $\tilde e$ orientation is $\A^1_{x'}$.

By the definition of contracted product, if we denote by $\msf a^*F$ the admissible $G$-cover
associated to $\msf a^*\tilde F$, then $\msf a^* F=\tau(b_F\cdot \msf a)\we^{\H^F}F$ as we wanted to prove.\fine

\begin{teo}\label{corocoro} 
Given a twisted $G$-cover $(\ssC,\phi)$ with 
associated admissible $G$-cover $F\to C$,
any ghost automorphism $\msf a\in \underline\Aut_C(\ssC)$ lifts to a ghost
automorphism of $(\ssC,\phi)$ if and only if the $1$-cochain
$b_F\cdot \msf a$ is in $\Ker\tau=\im\delta$ of sequence~(\ref{les_nor}).
\end{teo}
\proof After the proposition above, we have that $\phi\cong\msf a^*\phi$
if and only if $\tau(b_F\cdot \msf a)$ acts trivially via
the contracted product on $F$. We consider the restriction
$F_{\gen}\to C_{\gen}$ over the generic locus. We observe
that $F_{\gen}$ is an $\H^F$-torsor on $C_{\gen}$,
then we apply Proposition \ref{prop_fond} to obtain that
$\tau(b_F\cdot \msf a)\we^{\H^F}F_{\gen}=F_{\gen}$ if and only if
$b_F\cdot \msf a\in \im\delta^F$.
This is a necessary condition to have $\tau(b_F\cdot\msf a)\we^{\H^F}F=F$, but
it is also sufficient because $F_{\gen}$ completes uniquely to $F$.

It remains to prove that $\im\delta^F=\im\delta$. In particular we observe
that the contracted product does not act on the
$\delta$ morphism, so $\delta^F=\delta$ and the proof is concluded.\fine 

\begin{rmk}\label{rmk_z}
Given the dual graph $G$-cover $\tilde\Gamma\to\Gamma$ associated to $F\to C$,
and the contracted decorated graph $(\tilde \Gamma_0,r(-))$,
we recall the subcomplexes inclusion
$C^i(\tilde\Gamma_0;G)\subset C^i(\tilde\Gamma;G)$
 for $i=0,1$. We also consider
the exterior differential $\delta_0$ on $C^0(\tilde\Gamma;G)$, \cor{i.e.}~the restriction
of the $\delta$ operator to this group. Because of Proposition \ref{prop_imcap},
we have $\im(\delta_0)=C^1(\tilde\Gamma_0;G)\cap\im\delta$.
\end{rmk}

\begin{rmk}\label{rmk_prod}
Previously we obtained a characterization of the cochains
in $\im(\delta)$ that we could restate in our new setting. Indeed,
because of Proposition \ref{prop_imdelta.1}, an automorphism $\msf{a}\in S(\tilde \Gamma_0;r(-))$ is
an element of $\underline\Aut_C(\ssC,\phi)$ if and only if for every circuit $(\tilde e_1,\dots,\tilde e_k)$ in $\tilde\Gamma_0$
we have $\prod_{i=1}^k (b_F\cdot \msf a)(\tilde e_i)=1$.\newline
\end{rmk}

\subsection{Smooth points}
In Remark \ref{rmk_local} we discussed the fact that 
every point $[\ssC,\phi]\in \overline\R_{g,G}$ has a local picture isomorphic to
$\defo(\ssC)\slash\underline\Aut(\ssC,\phi)$.
This is a quotient of the form $\C^n\slash \G$ where $\G$
is a finite subgroup of $\GL(\C^n)$. In this setting we introduce
some automorphisms called quasireflections.

\begin{defin}[Quasireflection]\label{def_qr}
Any finite order complex automorphism $\msf h\in \GL(\C^n)$
is called a quasireflection if its fixed locus has dimension exactly $n-1$.
Equivalently, $\msf h$ is a quasireflection if, for an opportune choice of the basis,
we can diagonalize it as
$$\msf h=\diag(\xi,1,1,\dots,1),$$
where $\xi$ is a primitive root of the unit of order equal to the order of $\msf h$.
Given a finite group $\G\subset\GL(\C^n)$, we denote by $\QR(\G)$ the subgroup generated
by quasireflections.
\end{defin}

Quasireflections have the interesting property that any complex vector
space, quotiented by them, keeps being a smooth variety. In particular
if $\msf h\in GL(\C^n)$ is a quasireflection, the variety $\C^n\slash\msf h$ is isomorphic to $\C^n$.

\begin{prop}[see \cite{prill67}]\label{prop_qr}
Consider any vector space quotient
$V':=V\slash \G$,
where $V\cong \C^n$ is a complex vector space and $\G\subset \GL(V)$ is a finite group.
The variety $V'$ is smooth if and only if $\G$ is generated by quasireflections.
\end{prop}
Therefore, to find the smooth points of $\overline\R_{g,G}$,
by Proposition \ref{prop_qr}
we need to know when $\underline\Aut(\ssC,\phi)$ is generated by quasireflections.
We start by recalling the quasireflection analysis in the case of stable scheme theoretic curves.
\begin{defin}
Within a stable curve $C$, an elliptic tail is an irreducible component of geometric genus $1$ that
meets the rest of the curve in only one point called an elliptic tail node. Equivalently, $E$ is an elliptic tail if and only
if its algebraic genus is $1$ and $E\cap\overline{C\backslash E}=\{q\}.$
\end{defin}

An element $i\in\Aut(C)$ is an elliptic tail automorphism if there exists
an elliptic tail $E$ of $C$ such that $i$ fixes $E$ and his restriction to $\overline{C\backslash E}$ is the identity.
An elliptic tail automorphism of order $2$ is called an elliptic tail quasireflection (ETQR). In the literature ETQRs are
called elliptic tail involutions (or ETIs), we changed this convention in order to generalize
the notion.

\begin{rmk}\label{rmk_eti}
Every scheme theoretic curve of algebraic genus $1$ with one marked point has
exactly one involution $i$. Then there is a unique ETQR associated to every elliptic tail.

More precisely an elliptic tail $E$ could be
of two types. The first type is a smooth curve of geometric genus $1$
with one marked point, \cor{i.e.} an elliptic curve: in this case we have $E=\C\slash\Lambda$,
for $\Lambda$ integral lattice of rank $2$, the marked point is the origin,
and the only involution is the map induced by $x\mapsto -x$ on $\C$.
The second type is the rational line with one marked point and an autointersection
point: in this case we can write $E=\p^1\slash\{0\equiv \infty\}$, the marked point is the origin,
and the only involution is the map induced by $z\mapsto 1\slash z$ on $\p^1$.
\end{rmk}

From Remark \ref{rmk_def1} we have a coordinate system on $\defo(C)$
and on the canonical subscheme $\defo(C;\sing C)$.
Furthermore, the quotient of these two schemes has a splitting
$$\defo(C)\slash\defo(C;\sing C)\cong\bigoplus_{j=1}^{\delta}\A^1_{t_j}.$$
These coordinates systems on the space $\defo(C;\sing C)$ and $\defo(C)\slash\defo(C;\sing C)$
allow the detection of quasireflections. Indeed,
the diagonalizations of the $\msf a$-action on the two spaces determines a diagonalization of
the $\msf a$-action on the whole $\defo(C)$. Therefore, $\msf a$ is a quasireflection
if it acts non-trivially on exactly one coordinate of scheme $\defo(C;\sing C)$ or $\defo(C)\slash\defo(C;\sing C)$.
The following theorem by Harris and Mumford describes the action 
of the automorphism group $\Aut(C)$ on $\defo(C)$.

\begin{teo}[{See \cite[Theorem 2]{harmum82}}] \label{teo_eti}
Consider a stable curve $C$ of genus $g\geq 4$.
An element of $\Aut(C)$
acts as a quasireflection on $\defo(C)$ if and only if it is an ETQR. In particular, if $\eta\in \Aut(C)$ is an $ETQR$
acting non trivially on the tail $E$ with elliptic tail node $q_j$, then 
$\eta$ acts trivially on $\defo(C;\sing C)$, and its action on $\defo(C)\slash\defo(C;\sing C)$ is
$t_{j}\mapsto -t_{j}$ on the coordinate associated to $q_j$, and
 the identity $t\mapsto t$ on the remaining coordinates.\newline
\end{teo}

In Remark \ref{rmk_def2} we have seen that the deformations $\defo(C;\sing C)$ and
$\defo(\ssC;\sing \ssC)$ are canonically identified.
For the deformation of the nodes, the description is slightly different. We denote by $\delta$ the number
of nodes, by $r_1,\dots,r_{\delta}$
the order of the cyclic stabilizers in $\ssC$ of the nodes $q_1,\dots, q_{\delta}$ respectively.
Then,
$$\defo(\ssC)\slash\defo(\ssC;\sing \ssC)\cong\bigoplus_{j=1}^\delta \A^1_{\tilde t_j},$$
and every node comes with a flat representable morphism 
of Deligne-Mumford stacks, isomorphic to
$$
[\{x'y'=\tilde t_j\}\slash\mmu_{r_j}]\to \A^1_{\tilde t_j},$$
where the local stabilizer $\mmu_{r_j}$ acts by
$\xi\cdot (x',y',\tilde t_j)=(\xi x',\xi^{-1}y',\tilde t_j)$.
Also there exists a canonical morphism
$\A^1_{\tilde t_j}\to\A^1_{t_j}$ such that $(\tilde t_j)^{r_j}=t_j$.

\begin{rmk}\label{rmk_stacktail}
Consider a stack-theoretic curve $\msf E$ whose coarse space $E$
is a genus $1$ curve with a marked point. In the case of an elliptic
tail of a curve $\ssC$, the marked point is the point of intersection between $\msf E$
and $\overline{\ssC\backslash \msf E}$.

If $\msf E$ is an elliptic curve, then $\msf E=E$ and
the curve has exactly one involution $i_0$.
In case $\msf E$ is rational, its normalization
is the stack $\overline{\msf E}=[\p^1\slash\mmu_r]$, with $\mmu_r$ acting by multiplication,
and $\msf E=\overline{\msf E}\slash\{0\equiv \infty\}$. There exists a canonical involution
$\msf i_0$ in this case too: the pushforward of the inverse involution
on $\p^1$, \cor{i.e.}~$z\mapsto 1\slash z$.
We consider the autointersection node of $\msf E$ and its local picture $[\{x'y'=0\}\slash\mmu_r]$,
then the local picture of the same node in $E$ is $\{xy=0\}$ with $x=(x')^r$ and $y=(y')^r$.
Therefore the  $\msf i_0$-action is represented locally by $(x',y')\mapsto (y',x')$ and the product $x'y'$ is unchanged,
so $\msf i_0$ acts trivially on the smoothing coordinate $\tilde t$
associated to this node. We observe that of all the possible liftings of the canonical involution $i_0$ of $E$,
$\msf i_0$ is the only $\msf E$ involution acting trivially on $\tilde t$.

Given any twisted curve $\ssC$ with an elliptic tail $\msf E$ whose
elliptic tail node is called~$q$,
the construction above defines a canonical involution $\msf i_0$ on $\msf E$
up to non-trivial action on $q$.
\end{rmk}

\begin{defin}\label{def_etig}
An element $\msf i\in\underline\Aut(\ssC,\phi)$ is an ETQR if there exists
an elliptic tail $\msf E$ of $\ssC$ with elliptic tail node $q$, such that
the action of $\msf i$ on $\ssC\backslash\msf E$ is trivial, and the action
on $\msf E$, up to non-trivial action on $q$, is the canonical involution~$\msf i_0$.
\end{defin}

\begin{lemma}\label{lemma_qr}
Consider an element $\msf h$ of $\underline\Aut(\ssC,\phi)$. It acts as a quasireflection on $\defo(\ssC)$
if and only if one of the following is true:
\begin{enumerate}
\item the automorphism $\msf h$ is a ghost quasireflection,
\cor{i.e.}~an element of $\underline\Aut_C(\ssC,\phi)$ which moreover operates as
a quasireflection;
\item  the automorphism $\msf h$ is an ETQR, using the generalized Definition \ref{def_etig}.
\end{enumerate}
\end{lemma}
\proof 
We first prove the ``only if'' part. 
If $\msf h$ acts trivially on certain coordinates of $\defo(\ssC)$, \cor{a fortiori} we have that its coarsening $h$
acts trivially on the corresponding coordinates of $\defo(C)$. Therefore $h$ acts as the identity or as a quasireflection
on $\defo(C)$.
In the first case, $\msf h$ is a ghost automorphism and we are in case $(1)$. 
If $h$ acts as
a quasireflection, then it is a classical ETQR as we pointed out on Theorem \ref{teo_eti},
and it acts non-trivially on the coordinate associated to an elliptic tail node $q$.

As we know that the action of $h$ is trivial on $\defo(C;\sing C)$,
so is the action of $\msf h$.
It remains to know the action of $\msf h$ on the nodes with non-trivial stabilizer and other than $q$.
If the elliptic tail where $\msf h$ operates non trivially is a rational component with an autointersection
node $q_1$, by hypothesis $\msf h$ acts trivially on the universal deformation $\A^1_{\tilde t_1}$ of this node.
Therefore, the $\msf h$ restriction to the elliptic tail
has to be the canonical involution $\msf i_0$ (see Remark \ref{rmk_stacktail}).
For every node other than $q$ and $q_1$, if the local picture is $[\{x'y'=0\}\slash \mmu_{r_j}]$,
the action of $\msf h$ must be of the form
$$(x',y')\mapsto (\xi x',y')\equiv(x',\xi y')\ \ \m{ for some }\ \xi\in\mmu_{r_j}.$$
If $\xi\neq 1$ this gives a non-trivial action on the associated universal deformation $\A^1_{\tilde t_j}$,
against our hypothesis.
By Definition \ref{def_etig} this implies that $\msf h$ is an ETQR of $(\ssC,\phi)$.\newline

For the ``if'' part, we observe that a ghost quasireflection is automatically a quasireflection.
It remains to prove the case of point $(2)$. By definition of ETQR, its action on $\defo(\ssC)$
can be non-trivial
only on the components associated to the separating node $q$ of the tail. As a consequence
$\msf h$ acts as the identity or as a quasireflection.
The local coarse picture at $q$ is $\{xy=0\}$,
where $y=0$ is the branch lying on the elliptic tail. Then the action of $\msf h$ on the coarse space
is $(x,y)\mapsto(-x,y)$. Therefore the action is \cor{a fortiori} non trivial
on the coordinate associated to the stack node $q$ in $\defo(\ssC)$.\fine

\begin{defin}
For any stable curve $C$ we denote by $\QR(C)$ the subgroup of $\Aut(C)$
generated by classical ETQRs. For any twisted $G$-cover $(\ssC,\phi)$
we denote by $\QR(\ssC,\phi)$ the subgroup of $\underline\Aut(\ssC,\phi)$
generated by ETQRs, and by $\QR_C(\ssC,\phi)$ the subgroup of $\underline\Aut(\ssC,\phi)$
generated by ETQRs which moreover are ghosts.
\end{defin}

\begin{lemma}\label{keylemma}
Any element $h\in \QR(C)$ which
could be lifted to $\underline\Aut(\ssC,\phi)$, has a lifting in~$\QR(\ssC,\phi)$, too.
\end{lemma}
\proof By definition, $\underline\Aut(\ssC,\phi)$ is the set of automorphisms $\msf s\in\Aut(\ssC)$
such that $\msf s^*\phi\cong\phi$. Consider $h\in\QR(C)$ such that its decomposition
in ETQRs is $h=i_0i_1\cdots i_m$, and every $i_k$ acts non-trivially on an elliptic tail $E_k$.
Any lifting of $h$ is in the form $\msf h=\msf i_0\msf i_1\cdots \msf i_m\cdot \msf a$, where 
$\msf i_t$ is an ETQR acting non-trivially on a twisted elliptic tail $\msf E_k$, and $\msf a$ is a ghost
acting non-trivially only on nodes other than the elliptic tail nodes of the $\msf E_k$.
We observe that every $\msf i_k$ is a lifting
in $\Aut(\ssC)$ of~$i_k$. Moreover,
by construction, $\msf h^*\phi\cong\phi$ if and only if $\msf i_k^*\phi\cong \phi$ for
every $k$ and $\msf a^*\phi\cong \phi$. This implies that every $\msf i_k$ lies
in $\underline\Aut(\ssC,\phi)$, and therefore $\msf h\cdot\msf a^{-1}$ is a lifting of $h$
lying in $\QR(\ssC,\phi)$.\fine

We recall the short exact sequence (\ref{seq_aut}),
$$1\to\underline\Aut_C(\ssC,\phi)\to\underline\Aut(\ssC,\phi)\xrightarrow{\ \beta\ }\Aut'(C)\to 1$$
and introduce the group $\QR'(C)\subset\Aut'(C)$, generated by liftable quasireflections,
\cor{i.e.}~by those quasireflections $h\in\Aut(C)$ lying in $\im \beta$.
By Lemma \ref{keylemma}, $\QR'(C)=\beta(\QR(\ssC,\phi))$. Using also Lemma \ref{lemma_qr},
we obtain that the following is a short exact sequence
$$1\to \QR_C(\ssC,\phi)\to \QR(\ssC,\phi)\to \QR'(C)\to 1.$$

\begin{teo}\label{teoqr}
The group $\underline\Aut(\ssC,\phi)$ is generated by quasireflections
if and only if both $\underline\Aut_C(\ssC,\phi)$ and $\Aut'(C)$ are generated by
quasireflections.
\end{teo}
\proof
By combining the previous sequences,
$$1\to \underline\Aut_C(\ssC,\phi)\slash\QR_C(\ssC,\phi)\to
\underline\Aut(\ssC,\phi)\slash\QR(\ssC,\phi)\to \Aut'(C)\slash\QR'(C)\to 1.$$
The theorem follows.\fine

This gives a first important result for the moduli space of twisted $G$-covers $\overline\R_{g,G}$.
As we know that any point $[\ssC,\phi]\in\overline\R_{g,G}$ is smooth
if and only if the group $\underline\Aut(\ssC,\phi)$ is generated by quasireflections, then
the following theorem is straightforward.

\begin{teo}\label{teo_smooth}
Given a twisted $G$-cover $\phi\colon \ssC\to BG$
over a twisted curve $\ssC$ of genus $g\geq 4$ whose coarse space is $C$,
the point $[\ssC,\phi]$ of the moduli space $\overline\R_{g,G}$
is smooth if and only if the group $\Aut'(C)$ is generated by ETQRs
and the group of ghost automorphisms $\underline\Aut_C(\ssC,\phi)$ is
generated by quasireflections.
\end{teo}

We introduce two closed loci of $\overline\R_{g,G}$,
\begin{align*}
N_{g,G}&:=\left\{[\ssC,\phi]|\ \Aut'(C)\m{ is not generated by ETQRs}\right\},\\
H_{g,G}&:=\left\{[\ssC,\phi]|\ \underline\Aut_C(\ssC,\phi)\m{ is not generated by quasireflections}\right\}.
\end{align*}
We have by Theorem \ref{teo_smooth} that the singular locus $\sing \overline\R_{g,G}$ is
their union
\begin{equation}\label{eq_sing}
\sing\overline\R_{g,G}=N_{g,G}\cup H_{g,G}.
\end{equation}

\begin{rmk}
Consider the natural projection  $\pi\colon\overline\R_{g,G}\to\overline\M_g$,
then we have the inclusion
$N_{g,G}\subset\pi^{-1}\sing\overline\M_g$.
Indeed, we saw that $\QR'(C)=\Aut'(C)\cap \QR(C)$
and therefore $\Aut(C)=\QR(C)$ implies $\Aut'(C)=\QR'(C)$,
this means that $\left(\pi^{-1}\sing\overline\M_g\right)^c\subset (N_{g,G})^c$,
and taking the complementary we obtain the result.

We can interpret equality (\ref{eq_sing}) as the fact that the singular
locus is the union of two subloci: one coming from ``old'' singularities, the other
coming from data encoded only in the ghost structure of the twisted $G$-covers.\newline
\end{rmk}

The following lemma allows to characterize quasireflections in the ghost group.
Consider the decorated graph $(\Gamma(C),r(-))$ associated to a twisted $G$-cover $(\ssC,\phi)$,
and its contraction $(\Gamma_0,r(-))$ (see Definition \ref{gamma0}).

\begin{lemma}\label{ghost1}
Consider a ghost automorphism $\msf{a}$ in $\underline\Aut_C(\ssC)=S(\Gamma_0;r(-))$.
If $\msf a$ is a quasireflection in $\underline\Aut_C(\ssC,\phi)$
then
$\msf{a}(e)=1$ for all edges but one  that is a separating edge of~$\Gamma_0(\ssC)$.
\end{lemma}

\proof If $\msf{a}$ is a quasireflection in $\underline\Aut_C(\ssC,\phi)$, the value 
on all but one of the coordinates must be $0$.
Therefore $\msf{a}(e)=1\in\mmu_{r(e)}$ on all the edges but one, say $e_1$.
If there exists a preimage $\tilde e_1$ in $\E(\tilde\Gamma_0)$ that
is in any circuit $(\tilde e_1,\dots,\tilde e_k)$ of $\tilde \Gamma_0$ with $k\geq 1$,
then we have, by Remark \ref{rmk_prod},
that $\prod (b_F\cdot \msf{a})(\tilde e_i)=1$.
As $\msf{a}(e_1)\neq 1$, then $(b_F\cdot\msf a)(\tilde e_1)\neq 1$ and therefore
 there exists $i>1$ such that $(b_F\cdot\msf{a})(\tilde e_i)\neq 1$ too.
 This would imply that, if $e_i$ is the image in $\Gamma_0$ of $\tilde e_i$,
 then $\msf a(e_i)\neq 1$, contradiction.
Thus $\tilde e_1$ is not in any circuit, then it is a separating edge
and so is $e_1$.

Reciprocally,
consider an automorphism $\msf{a}\in S(\Gamma_0;r(-))$ such
that there exists an oriented separating edge $e_1$ with the property that
$\msf{a}(e)= 1 $ for every $e$ in $\E\backslash\{e_1, \bar e_1\}$ and
$\msf{a}(e_1)$ is a non-zero element of $\mmu_{r(e_1)}$. Then for every circuit $(\tilde e_1',\dots,\tilde e_k')$
of $\E(\tilde \Gamma_0)$, we have
$\prod (b_F\cdot\msf{a})(\tilde e_i')=1$ and so $\msf{a}$ is in $\underline\Aut_C(\ssC,\phi)$
by Theorem \ref{corocoro}.\fine

\section{Non-canonical singularities}\label{ncs}

\subsection{Characterization of the non-canonical locus}
In order to detect the singularity canonicity, we need a tool called age invariant. 
After its introduction we will be able to prove the bipartition of $\sing^{\nc}\overline\R_{g,G}$.

\subsubsection{The age invariant}
Consider the case of a vector space quotient $V\slash\G$.
In the case of the group $\G$ not being generated by quasireflections, we need
another tool to distinguish between canonical singularities and non-canonical singularities.
The age is a positive function $\G\to \Q$.

\begin{defin}[Age]\label{defage}
Consider a $\G$-representation $\rho\colon \G\to \GL(V)$.
For any element $\msf h\in \G$ of order $r$, there exists a diagonalization
$\msf h=\diag(\xi_r^{a_1},\xi_r^{a_2},\dots,\xi_r^{a_n})$, where $\xi_r=\exp(2\pi i\slash r)$
is a privileged $r$th root of the unit and $0\leq a_i<r$ for any $i=1,\dots n$. In this setting
$$\age(\msf h)=\frac{1}{r}\sum_{i=1}^na_i.$$
\end{defin}

\begin{defin}[Junior group]\label{def_junior}
A finite group $\G\subset \GL(\C^m)$ that contains no quasireflections
 is called junior if the image of the age function intersects the open
interval~$]0,1[$,
$$\age \G\cap\ ]0,1[\ \neq \varnothing.$$
The group $\G$ is called senior if the intersection is empty.
\end{defin}

\begin{rmk}\label{rmkchoice}
The definition of age depends on the non-canonical choice of a privileged
root $\xi_r$, but the image $\age(\G)\subset \Q$ does not depend on this choice.
Therefore junior and senior group are well defined.
\end{rmk}

\begin{prop}[Age criterion, see	 \cite{reid80}]\label{agec}
Consider any vector space quotient
$V':=V\slash \G$,
where $V\cong \C^n$ is a complex vector space and $\G\subset \GL(V)$
is a finite group containing no quasireflections. Then $V'$
has a non-canonical singularity if and only if $\G$ is junior.
\end{prop}

We will use the Age Criterion to find
non-canonical singularities by the study of group $\underline\Aut(\ssC,\phi)$ action
on $\defo(\ssC,\phi)$.
We point out that to satisfy the hypothesis of Age Criterion, it is necessary
for $\underline\Aut(\ssC,\phi)$ to be quasireflection free. As this is often not the case,
the following lemma is necessary to represent the same singularity by
a group with no quasireflections.

\begin{prop}[see \cite{prill67}]\label{propt}
Consider a finite subgroup $\G\subset\GL(\C^n)$. There exists an isomorphism 
$u \colon \C^n\slash \QR(\G)\to\C^n$
and a finite subgroup $\mfk K\subset\GL(\C^n)$ isomorphic to the quotient $\G\slash\QR(\G)$, such that the following
diagram is commutative.
\[
\begin{CD}
\C^n@>>> \C^n\slash \QR(\G)@>u>>\C^n\\
@VVV @VVV @VVV\\
\C^n\slash \G@>\cong>> (\C^n\slash \QR(\G))\slash (\G\slash\QR(\G))@>\cong>> \C^n\slash {\mfk K}\\
\end{CD}
\]\newline
\end{prop}

\subsubsection{$T$-curves and $J$-curves}
We introduce two closed loci which are central
in our description.
\begin{defin}[$T$-curve]
A twisted $G$-cover $(\ssC,\phi)$ is a $T$-curve if
there exists
an automorphism $\msf a\in \underline\Aut(\ssC,\phi)$ such that its 
coarsening $a$ is an elliptic tail automorphism of order~$6$.
The locus of $T$-curves in $\overline\R_{g,G}$ is denoted by $T_{g,G}$.
\end{defin}

\begin{defin}[$J$-curve]
A twisted $G$-cover $(\ssC,\phi)$ is a $J$-curve if
the group 
$$\underline\Aut_C(\ssC,\phi)\slash\QR_C(\ssC,\phi),$$
which is the group of ghosts quotiented by its subgroup of quasireflections, is junior.
The locus of $J$-curve in $\overline\R_{g,G}$ is denoted by $J_{g,G}$.
\end{defin}

\begin{teo}\label{teo_jt}
For $g\geq 4$, the non-canonical locus of $\overline\R_{g,G}$ is the union
$$\sing^{\nc}\overline\R_{g,G}=T_{g,G}\cup J_{g,G}.$$
\end{teo}

\begin{rmk}
We observe that \cite[Theorem 2.44]{chiofar12}, affirms
exactly that in the case $G=\mmu_\ell$ with $\ell\leq 6$ and $\ell\neq 5$, the $J$-locus
$J_{g,\mmu_\ell}$ is empty for every genus $g$, and therefore $\sing^{\nc}\overline\R_{g,\mmu_\ell}$
coincides with the $T$-locus for these values of $\ell$.
\end{rmk}


We introduce the notion of $\star$-smoothing, following
\cite{harmum82} and~\cite{lud10}. 
\begin{defin}
Consider a twisted $G$-cover $(\ssC,\phi)$
and  a junior  automorphism
 $\msf{a}\in\underline\Aut(\ssC,\phi)\slash\QR(\ssC,\phi)$,
we say that  the triple $(\ssC,\phi,\msf a)$ is $\star$-smoothable if 
\begin{itemize}
\item on the coarse curve $C$ there exists a cycle of $m$ non-separating nodes $q_0,\dots,q_{m-1}$,
\cor{i.e.}~we have $\msf a(q_i)=q_{i+1}$ for all $i=0,1\dots,m-2$ and $\msf{a} (q_{m-1})=q_0$;
\item the action of $\msf a^m$ over
the coordinate associated to every node is trivial. Equivalently, $\msf a^m(\tilde t_{q_i})=\tilde t_{q_i}$
for all $i=0,1,\dots,m-2$, where $\tilde t_{q_i}$ is the coordinate on $\defo(\ssC,\phi)$ associated
to the $q_i$-smoothing (see Remark \ref{rmk_def2}).
\end{itemize}
\end{defin}

If $(\ssC,\phi,\msf a)$ is $\star$-smoothable, there exists
a deformation $(\ssC',\phi',\msf a')$
that smooths the $m$ nodes and
with $\msf a'\in\underline\Aut(\ssC',\phi')$.
Moreover, this deformation preserves the age of the $\msf a $-action on $\defo(\ssC,\phi)\slash\QR$. Indeed,
the eigenvalues of $\msf a$ are a discrete and locally constant set, thus constant by
deformation.
The $T$-locus and the $J$-locus are closed by $\star$-smoothing,
\cor{i.e.}~if the deformation $(\ssC',\phi')$ above is a $T$-curve or a $J$-curve,
then $(\ssC,\phi)$ is a $T$-curve or a $J$-curve.

Therefore in proving Theorem \ref{teo_jt}, we can suppose that every triple $(\ssC,\phi,\msf a)$
 that we consider is $\star$-\cor{rigid}, \cor{i.e.}~non-$\star$-smoothable.
Indeed, if there exists
a junior $\star$-smoothable automorphism $\msf a\in \underline\Aut(\ssC,\phi)$, 
we smooth it until we obtain a rigid triple $(\ssC',\phi',\msf a')$ of the same age.
Then, if $(\ssC',\phi')$ is a $T$-curve or a $J$-curve, the same is true for $(\ssC,\phi)$.\newline

\proof[Proof of Theorem \ref{teo_jt}]
We will show in eight steps that if the group
$\underline\Aut(\ssC,\phi)\slash\QR\left(\ssC,\phi\right)$
is junior, and $(\ssC,\phi)$ is not a $J$-curve, then it is a $T$-curve.
After the Age Criterion \ref{agec} and Proposition \ref{propt}, this will prove
Theorem \ref{teo_jt}. 
From now on we work under the hypothesis that $\msf a\in\underline\Aut_C(\ssC,\phi)\slash\QR$
is a non-trivial automorphism aged less than $1$, 
that $(\ssC,\phi)$ is not a $J$-curve and $(\ssC,\phi,\msf a)$ is $\star$-rigid.

In steps $1$ and $2$ we fix the setting and prove two useful lemmata. In step $3$
we prove that all the nodes of $\ssC$ are fixed by $\msf a$ except at most $2$ of them
which are exchanged.
In step $4$ we show that every irreducible component $\msf Z\subset \ssC$ is fixed by $\msf a$.
In step $5$ we can therefore conclude that there are no couple of exchanging nodes.
In step $6$ and $7$ we study the action of $\msf a$ on the irreducible components
of $\ssC$ and the contributions to $\age \msf a$. Finally we prove the result in step $8$.\newline

\cor{Step 1.}
Consider the contracted decorated graph $(\Gamma_0,r(-))$ of $(\ssC,\phi)$.
As before, we call $E_{\sep}$ the set of separating edges of $\Gamma_0$.
As stated in Remark \ref{rmk_def2},
we have the following splitting,
\begin{equation}\label{eq_split}
\Def(\ssC,\phi)\cong\defo(\ssC;\sing\ssC)\oplus\bigoplus_{e\in E_{\sep}}\A_{\tilde t_e}\oplus \bigoplus_{e'\in E\backslash E_{\sep}}\A_{\tilde t_{e'}},
\end{equation}
where $\tilde t_e$ is a coordinate parametrizing the smoothing of the node associated to the edge $e$.
In particular for every vector subspace $V\subset \Def(\ssC,\phi)$
and every automorphism $\msf a$ of $(\ssC,\phi)$, 
we denote by $\age(\msf a|V)$ the
age of the restriction $\msf a|_V$. If $\msf Z$ is a subcurve of $\ssC$, then there exists
a canonical inclusion $\Def(\msf Z)\subset \Def(\ssC)$, and we define
$\age(\msf a|\msf Z):=\age(\msf a|\Def(\msf Z))$.\newline

Every ghost automorphism in $\underline\Aut(\ssC,\phi)$ fixes
the three summands of (\ref{eq_split}).
Moreover, every quasireflection acts only on the summand $\bigoplus_{e\in E_{\sep}}\A_{\tilde t_e}$ by Lemmata
 \ref{lemma_qr} and \ref{ghost1}.
As a consequence,
by Propostion~\ref{propt}, the group $\underline\Aut(\ssC,\phi)\slash\QR$ acts on
\begin{equation}\label{eq_space}
\frac{\Def(\ssC,\phi)\slash\QR}{\Def(\ssC;\sing\ssC)}\cong\left( \frac{\bigoplus_{e\in E_{\sep}}\A_{\tilde t_e}}{\QR(\ssC,\phi)}\right)\oplus
 \bigoplus_{e'\in E\backslash E_{\sep}}\A_{\tilde t_{e'}}.
\end{equation}

Every quasireflection acts on exactly one coordinate $\tilde t_e$
with $e\in E_{\sep}$. We rescale all the coordinates $\tilde t_e$
by the action of $\QR(\ssC,\phi)$. We call $\tau_e$,
for $e\in E(\Gamma_0)$, the new set of coordinates.
Obviously $\tau_{e'}=\tilde t_{e'}$ if $e'\in E(\Gamma_0)\backslash E_{\sep}$.\newline

\cor{Step 2.} We show two lemmata about the age contribution of
the $\msf a$-action on nodes, that we call \emph{aging} on nodes.

\begin{defin}[coarsening order]
If $\msf a\in\underline\Aut(\ssC,\phi)$ and $a$ is its coarsening, then we define 
$$\cord\msf a:=\ord a.$$
\end{defin}
The coarsening order is the least integer $n$ for which $\msf a^n$ is a ghost automorphism.

\begin{lemma}\label{lemma_bound}
Suppose that $\msf Z\subset \ssC$ is a subcurve of $\ssC$ such that $\msf a(\msf Z)=\msf Z$ and 
$q_0,\dots q_{m-1}$ is a cycle, by $\msf a $, of nodes in $\msf Z$. Then
we have the following inequalities:
\begin{enumerate}
\item $\age(\msf a|{\msf Z})\geq \frac{m-1}{2}$;
\item  if the nodes $q_0,\dots,q_{m-1}$ are non-separating, $\age (\msf a)\geq \frac{m}{\ord\left(\msf a|{\msf Z}\right)}+\frac{m-1}{2}$;
\item if $\msf a^{\cord \msf a}$ is a senior ghost, we have $\age (\msf a)\geq \frac{1}{\cord(\msf{a})}+\frac{m-1}{2}$.
\end{enumerate}
\end{lemma}
\proof
We call $\tau_0,\tau_1,\dots,\t_{m-1}$ the coordinates associated to nodes
$q_0,\dots,q_{m-1}$ respectively. By hypothesis, $\msf a ( \tau_0)=c_1\cdot  \tau_1$ and
$\msf a^i(\t_0)=c_i\cdot\t_i$ for all $i=2,\dots,m-1$,
where the $c_i$ are complex numbers. If $n'=\ord(\msf a|_{\msf Z})$,
we have 
$$\msf a^m(\t_0)=\xi_{n'}^{u m}\cdot \t_0$$
 where $\xi_{n'}$
is a primitive $n'$th root of the unit and $u$ is an integer such that 
$0 \leq u <n'\slash m$.
The integer $u$ is called \cor{exponent} of the cycle $(q_0,\dots,q_{m-1})$
with respect to the curve~$\msf Z$.
Observe that $\msf a(\t_{i-1})=(c_i\slash c_{i-1})\cdot \t_i$ and $\msf a^m (\t_i)=\xi_{n'}^{u m}\cdot\t_i$
for every $i$.

We can explicitly write the eigenvectors for the action of $\msf a$ on the coordinates
$\t_0,\dots,\t_{m-1}$. Set $d:={n'}\slash m$ and $ b:=sd+u$ with $0\leq s<m$, and consider the vector
$$\mfk v_b:=(\t_0=1,\ \t_1=c_1\cdot \xi_{n'}^{-b},\dots,\ \t_i=c_i\cdot\xi_{n'}^{-ib},\dots).$$
Then $\msf a (\mfk v_b)=\xi_{n'}^b\cdot \mfk v_b$. The contribution to the age of the eigenvalue $\xi_{n'}^b$ is
$b\slash n'$, thus we have
$$\age\msf a \geq \sum_{s=0}^{m-1}\frac{sd+u}{n'}=\frac{mu}{n'}+\frac{m-1}{2},$$
proving point $(1)$.

If the nodes are non-separating, as we are supposing that $(\ssC,\phi,\msf a)$ is $\star$-rigid,
we have $u\geq 1$ and the point $(2)$ is proved.

Suppose that $\msf a$ has order $n=\ord \msf a$ and its action on $\ssC$ has $j$ nodes cycles
of order $m_1, m_2,\dots, m_j$ and exponents respectively 
$u_1,\dots, u_j$ with respect to $\ssC$. If $k=\cord \msf a$, then
 $\msf a^k$ fixes every node, then we consider the coordinate $\t_i$ of a node
of the first cycle and we have
$$\msf a^k(\t_i)=\xi_n^{w\cdot k}\cdot\t_i,$$
where $w$ is an integer such that $0\leq w <n\slash k$. 
Repeating the same operation for every cycle we obtain
another series of integers $w_1,w_2,\dots, w_j$.
Therefore the age of $\msf a^k$ is 
$$\age \msf a^k=\sum_{i=1}^j\frac{m_iw_ik}{n},$$ and it is greater
or equal to $1$ by hypothesis.

We observe that $m_i$ divides $k$ for all $i=1,\dots, j$, and
$$u_i\cdot m_i\cdot \frac{k}{ m_i}\equiv w_i\cdot k\mod n.$$
This implies that $u_i\geq w_i$ for every $i$.

By the point $(2)$, the age of $\msf a$ on the $i$th cycle
is bounded from below by 
$m_iu_i\slash n+(m_i-1)\slash 2$.
As a consequence
$$\age \msf a\geq \sum_{i=1}^j\left(\frac{m_iu_i}{n}+\frac{m_i-1}{2}\right)\geq \sum_{i=1}^j\left(\frac{m_iw_i}{n}+\frac{m_i-1}{2}\right) \geq \frac{1}{k}+\frac{m_1-1}{2}.$$ \fine

\cor{Step 3.} Because of  Lemma \ref{lemma_bound},
if the automorphism $\msf a$ induces a cycle of $m$ nodes,
then this cycles contributes by at least $\frac{m-1}{2}$ to the aging of $\msf a$.
Therefore, as $\msf a$ is junior,
all the nodes of $\ssC$ are fixed except at most two of
them, that are exchanged. Moreover, if a pair of non-fixed
nodes exists, they contribute by at least $1\slash 2$.\newline

\cor{Step 4.} Consider an irreducible component $\msf Z\subset \ssC$, we want to prove $\msf a(\msf Z)=\msf Z$.
Suppose there is a cycle of irreducible components $\msf C_1, \dots,\msf C_m$  with $m\geq 2$ such that
$\msf a(\msf C_i)=\msf C_{i+1}$ for $i=1,\dots,m-1$, and $\msf a(\msf C_m)=\msf C_1$.
We call $\overline\ssC_i$ the normalizations of these components, and
$D_i$ the preimages of $\ssC$ nodes on $\overline\ssC_i$.
We point out that this construction implies that $(\overline\ssC_i,D_i)\cong (\overline\ssC_j, D_j)$ for all $i,j$.
Then,
an argument of \cite[p.34]{harmum82} shows that the action of $\msf a$ on
$\defo(\ssC;\sing \ssC)$
gives a contribution of at least $k\cdot(m-1)\slash2$ to $\age\msf a$, where
$$k=\dim H^1(\overline\ssC_i,T_{\overline\ssC_i}(-D_i))=3g_i-3+|D_i|.$$

This gives us two cases for which $m$ could be greater than $1$
with still a junior age: $k=1$ and $m=2$ or $k=0$.

If $k=1$ and $m=2$, we have $g_i=0$ or $1$ for $i=1,2$.
Moreover, the aging of at least $1\slash 2$ sums to
another aging of $1\slash 2$ if there is a pair of non-fixed nodes.
As $\msf a$ is junior, we conclude that $\ssC=\ssC_1\cup\msf a(\ssC_1)$
but this implies $g(\ssC)\leq 3$, contradiction.

If $k=0$, we have $g_i=1$ or $g_i=0$, the first is excluded because
it implies $|D_i|=0$ but the component must intersect the curve somewhere.
Thus, for every component in the cycle, the normalization $\overline\ssC_i$ is the
projective line $\p^1$ with $3$ marked points. We have two cases:
the component $\ssC_i$ intersects $\overline{\ssC\backslash\ssC_i}$ in $3$ points
or in $1$ point, in the second case $\ssC_i$ has an autointersection node
and $\ssC=\ssC_1\cup\msf a(\ssC_1)$, which is a contradiction because $g(\ssC)< 4$.
It remains the case in the image below.

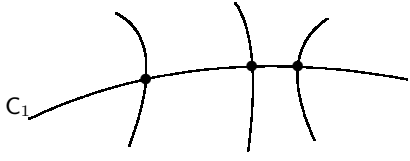
\begin{figure}[h]
 \begin{picture}(400,50)(40,0)
   \qbezier(187,10)(250,40)(330,25)
   \qbezier(220,50)(240,40)(225,0)
   \qbezier(265,54)(275,40)(270,-2)
   \qbezier(300,48)(280,35)(295,2)
   
   \put(178,12){\m{\footnotesize{$ \ssC_1$}}}
   \put(231.2,25.3){\circle*{4}}
   \put(271.3,30){\circle*{4}}
   \put(288.8,30){\circle*{4}}
 
 \end{picture}
 \caption{Case with $ \ssC_1\cong \p^1$ and $3$ marked points}
\end{figure} 

As $\ssC_1,\ssC_2,\dots,\ssC_m$ are moved by $\msf a$, every
node on $\ssC_1$ is transposed with another one or is fixed 
with its branches interchanged.
If at least two nodes are transposed we have an age contribution bigger or equal to $1$
by Lemma \ref{lemma_bound}.
If only one node is transposed we have two cases. In the first case
$\ssC=\ssC_1\cup \msf a(\ssC_1)\cup\ssC_2\cup\msf a(\ssC_2)$,
where $\ssC_2$ intersects only the component $\ssC_1$ and in exactly one point.
If $g(\ssC_2)\geq 2$, then the age is bigger than $1$, if $g(\ssC_2)<2$, then
$g(\ssC)\leq 3$, contradiction.

In the second case, $\ssC=\ssC_1\cup\msf a(\ssC_1)\cup \ssC_2$ where
$\ssC_2$ intersects $\ssC_1$ and $\msf a(\ssC_1)$, both in exactly one point.
If $g(\ssC_2)<2$ we have another genus contradiction. By the results
of \cite[p.28]{harmum82}, to have $\age(\msf a | \ssC_2)<1$ we must have
$g(\ssC_2)=2$ and the coarsening of $\msf a $ has order $2$. 
Therefore by Lemma \ref{lemma_bound} point (3), $\msf a$ has
age bigger or equal to $1$.\newline

\cor{Step 5.} We prove that every node is fixed by $\msf a$.
Consider the normalization $\msf{nor}\colon\bigsqcup_i\overline\ssC_i\to\ssC$
already introduced.
If the age of $\msf a$ is lower than $1$, \emph{a fortiori} 
we have $\age (\msf a|{\overline\ssC_i})<1$ for all~$i$.
In \cite[p.28]{harmum82}
there is a list of those smooth stable curves for which there exists
a non-trivial junior action.

\begin{enumerate}[label=\roman{*}., ref=(\roman{*})]
\item The projective line $\p^1$ with $\msf a\colon z\mapsto (-z)$ or $(\xi_4 z)$;
\item an elliptic curve with $\msf a$ of order $2,\ 3,\ 4$ or $6$;
\item an hyperelliptic curve of genus $2$ or $3$ with $\msf a$ the hyperelliptic involution;
\item a bielliptic curve of genus $2$ with $\msf a$ the canonical involution.
\end{enumerate}
We observe that the order of the $\msf a$-action on these components is always $2,\ 3,\ 4$ or~$6$.
As a consequence, if $\msf a$ is junior, then $n=\cord\msf a=2,3,4,6$ or~$12$,
as it is the greatest common divisor between the $\cord\left(\msf a|{\overline\ssC_i}\right)$.

First we suppose $\ord\msf a>\cord\msf a$, thus $\msf a^{\cord\msf a}$ is a ghost and it must
be senior. Indeed, if $\msf a^{\cord\msf a}$ is aged less than $1$,
then $(\ssC,\phi)$ admits junior ghosts,
contradicting our assumption. By point $(3)$ of
Lemma \ref{lemma_bound}, if there
exists a pair of non-fixed nodes, we obtain an aging of $1\slash n+1\slash2$ on
node coordinates.
If  $\ord \msf a=\cord\msf a$ the bound is even greater.
As every component is fixed by $\msf a$, the two nodes are non-separating, and by point $(2)$
of Lemma \ref{lemma_bound} we obtain an aging of $2\slash n+1\slash2$.

If $\overline\ssC_i$ admits an automorphism of order $3,\ 4$ or $6$, 
by a previous 
analysis of Harris and Mumford (see \cite{harmum82} again),
this yields
an aging  of, respectively, $1\slash 3$,\ $1\slash 2$ and $1\slash 3$ on $H^1(\overline\ssC_i,T_{\overline\ssC_i}(-D_i))$.

These results combined, show that a non-fixed pair of nodes gives
an age greater than $1$. Thus, if $\msf a$ is junior,
every node is fixed.\newline

\cor{Step 6.} We study the action of $\msf a$ separately on every irreducible
component. The $\msf a$-action is non-trivial on at least one component~$\ssC_i$, and this component must lie
in the list above. 

In case (i), $\overline\ssC_i$ has at least $3$ marked points
because of the stability condition. Actions of type $x\mapsto \xi x$
have two fixed points on $\p^1$, thus
at least one of the marked points is non-fixed.
A non-fixed preimage of a node has order $2$, thus
 the coarsening $a$ of $\msf a$ is the involution $z\mapsto -z$. Moreover,
$\ssC_i$ is the autointersection of the projective line and
$\msf a$ exchanges the branches of the node. 
Therefore $\msf a^2|_{\overline\ssC_i}$ is a ghost automorphism
of $\overline\ssC_i$. As a direct consequence of Theorem \ref{corocoro} and Remark \ref{rmk_prod},
the action of $\msf a^2$ on the coordinate associated to
the autointersection node, is trivial.
Therefore the action of $\msf a^2$ on
the same coordinate gives an aging of $0$ or $1\slash 2$, by
$\star$-rigidity it is $1\slash 2$.

The analysis for cases (iii) and (iv) is identical to that
developed in \cite{harmum82}: the only possibility of a junior action is the case
of an hyperelliptic curve $\msf E$ of genus $2$ intersecting $\overline{\ssC\backslash\msf E} $
in exactly one point, whose hyperinvolution gives an 
aging of $1\slash 2$ on $H^1(\overline\ssC_i,T_{\overline\ssC_i}(-D_i))$.

Finally, in case (ii), we use again the analysis of \cite{harmum82}.
The elliptic component $\msf E$ has
$1$ or $2$ point of intersection with $\overline{\ssC\backslash\msf E}$.
If there is $1$ point of intersection, elliptic tail case,
for a good choice of coordinates
the coarsening $a$ acts as $z\mapsto \xi_nz$,
where $n$ is $2,3,4$ or $6$. The aging is, respectively,
$0,\ 1\slash3,\ 1\slash2,\ 1\slash3.$ If there are $2$ points of
intersection, elliptic ladder case, the order of $\msf a$ on $\msf E$ must be $2$ or $4$ and the
aging respectively $1\slash 2$ or $3\slash 4$.\newline

\cor{Step 7.} 
Resuming what we saw until now, if $\msf a$ is a junior
automorphism of $(\ssC,\phi)$, $a$ its coarsening
and $C_1$ an irreducible component of $C$,
then we have one of the following:
\begin{enumerate}
\item [A.] component $C_1$ is an hyperelliptic tail, crossing the curve in one point,
with $a$ acting as the hyperelliptic involution and aging $1\slash 2$ on $H^1(\overline C_1,T_{\overline\ssC_i}(-D_1))$;
\item [B.] component $C_1$ is a projective line $\p^1$ autointersecting itself, crossing the curve in one point,
with $a$ the involution which fixes the nodes, and aging $1\slash 2$;
\item [C.] component $C_1$ is an elliptic ladder, crossing the curve in two points,
with $a$ of order $2$ or $4$ and aging respectively $1\slash2$ or $3\slash 4$;
\item [D.] component $C_1$ is an elliptic tail, crossing the curve
in one point, with $a$ of order $2,3,4$ or $6$ and aging $0,1\slash3,1\slash 2$ or $1\slash 3$;
\item [E.] automorphism $a$ acts trivially on $C_1$ with no aging.\newline
\end{enumerate}

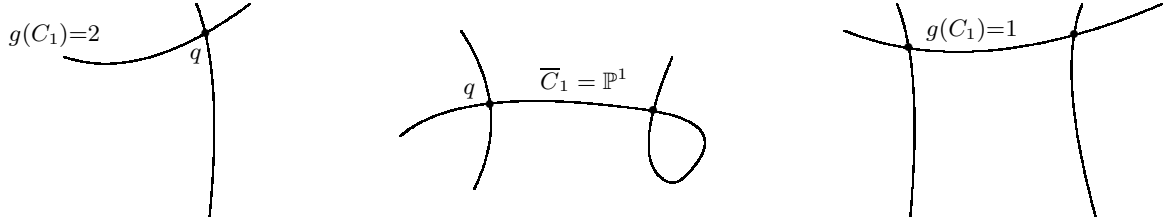
\begin{figure}[h]
 \begin{picture}(400,70)(40,0)
 \put(20,-10){
 \qbezier(20,70)(50,60)(90,90)
 \qbezier(70,90)(80,65)(75,10)
 
 \put(73.4,79.3){\circle*{3}}
 
 \put(-1,76){\m{\footnotesize{$g(C_1)$$=$$2$}}}
 
 \put(67.5,69.4){\m{\footnotesize{$q$}}}}
 
 \put(-20,0){
 \qbezier(187,30)(210,50)(280,40)
 \qbezier(280,40)(315,35)(295,15)
 \qbezier(295,15)(290,10)(285,15)
 \qbezier(285,15)(275,25)(290,60)
 
 \qbezier(210,70)(230,40)(215,10)
 
 \put(221,42.4){\circle*{3}}
 \put(211,46){\m{\footnotesize{$q$}}}
 
 \put(282.8,39.9){\circle*{3}}

\put(240,48){\m{\footnotesize{$\overline C_1=\p^1$}}} }

 \put(335,-10){
 \qbezier(0,80)(50,60)(120,90)
 \qbezier(20,90)(30,65)(25,10)
 \qbezier(90,90)(80,65)(95,10)
 
 \put(87,78.6){\circle*{3}}
 \put(24.4,74){\circle*{3}}
 
 \put(31,78){\m{\footnotesize{$g(C_1)$$=$$1$}}}}
 
 \end{picture}
 \caption{Components of type A, B and C.}
\end{figure} 

We rule out cases (A), (B) and (C).
At first we suppose
there is a component of type (A) or (B). For genus reasons, the component
intersected in both cases must be of type (E).
We study the local action on
the separating node $q$. The local picture at $q$ is 
$\left[\{x'y'=0\}\slash\mmu_r\right]$.
The smoothing
of the node is given by the stack $\left[\{xy=\tilde t_q\}\slash\mmu_r\right]$.
Consider the action of the automorphism $\msf a$ at the node, as the coarsening
of $\msf a$ has order $2$, then $\msf a\colon\tilde t_q\mapsto \varsigma\cdot \tilde t_q$
and $\varsigma^2\in\mmu_r$. Therefore $\msf a^2$ acts as the identity or as a quasireflection
of factor $\varsigma^2$. 
Thus $\t_{q}=\tilde t_{q}^{r'}$ where $r'|r$ is the order of $\varsigma^2$.
Therefore the action of
$\msf a$ on $\A^1_{\tilde t_q}\slash\QR=\A^1_{\t_q}$ is $\t_{q}\mapsto \varsigma^{r'}\cdot\t_{q}=-\t_q$. The
additional age contribution is $1\slash 2$, ruling out this case.

In case there is a component of type (C), if its nodes
are separating, then one of them must intersect a component
of type (E) and we use the previous idea. In case nodes are non-separating, 
we use Lemma \ref{lemma_bound}. If $\ord \msf a>\cord \msf a$, then $\msf a^{\cord\msf a}$ is
a senior ghost because $(\ssC,\phi)$ is not a $J$-curve, thus by
point $(3)$ of the lemma there is an aging of $(1\slash\cord\msf a)$ on the node coordinates.
If $\ord\msf a=\cord\msf a$, the bound is even greater, as by point $(2)$
we have an aging of $(2\slash\cord\msf a)$. We observe that $\cord\msf a=2,4$ or $6$,
and in case $\cord\msf a=6$ there must be a component of type (E).
Using additional contributions listed above we rule out
the case (C).\newline

\cor{Step 8.} We proved that $\ssC$ contains
components of type (D) or (E), \cor{i.e.}~the automorphism
$\msf a$ acts non-trivially only on elliptic tails.
If $q$ is the elliptic tail node, there are
two quasireflections acting on the coordinate $\tilde t_{q}$:
a ghost automorphism associated to this node
and the elliptic tail quasireflection. If the order of the local stabilizer is $r$,
then $\t_{q}=\tilde t_{q}^{2r}$.

If $\ord \msf a=2$ we are in the ETQR case, this action is a quasireflection
and it contributes to rescaling the coordinate $\tilde t_{q}$.

If $\ord\msf a=4$, the action on the (coarse) elliptic tail is $z\mapsto \xi_4z$. The space
$H^1(\overline C_i,T_{\overline C_i}(-D_i))$ is the space of $2$-forms
$H^0(\overline C_i,\omega^{\xx 2}_{\overline C_i})$: this space is generated by $dz^{\xx 2}$
and the action of $\msf a$ is $dz^{\xx 2}\mapsto\xi_4^2dz^{\xx2}$.
Moreover, if the local picture of the elliptic tail node is $[\{x'y'=0\}\slash\mmu_r]$,
then $\msf a\colon(x',y')\mapsto(\zeta x',\varrho y')$ such that
$\zeta^r=\xi_4$ and $\varrho^r=1$. As a consequence $\msf a\colon\tilde t_{q}\mapsto \zeta\cdot\varrho\cdot\tilde t_{q}$
and therefore $\t_{q}\mapsto \xi_2\t_{q}$.
Then, $\age \msf a=1\slash2+1\slash2$, proving the seniority of $\msf a$.

If $\msf E$ admits an automorphism $\msf a$ of order $6$,
 the action
on the (coarse) elliptic tail is $\msf a\colon z\mapsto \xi_6^kz$.
Then $dz^{\xx 2}\mapsto \xi_3^kdz^{\xx2} $ and $\t_{q}\mapsto\xi_3^k\t_{q}$.
For $k=1,4$ we have age lower than~$1$.\newline

If $(\ssC,\phi)$ is not a $J$-curve, 
we have shown that
 the only case where an automorphism $\msf a$
 in $\underline\Aut(\ssC,\phi)\slash\QR$ is junior,
is when its coarsening $a$ is an elliptic tail automorphism of order~$6$.\fine

\subsection{The $J$-locus in the case $S_3$}
We consider the case of $J_{g,S_3}$ and prove, thanks to the tools we developed,
that this locus is empty, that is the following.

\begin{teo}\label{teo_ncs3}
If $G$ is the symmetric group
$S_3$, then the non-canonical locus coincides
with the $T$-locus,
$$\sing^{\nc}\overline\R_{g,S_3}=T_{g,S_3}.$$
In particular, a point $[\ssC,\phi]$ is a non-canonical singular point if and only if
there exists an automorphism $\msf a\in\underline\Aut(\ssC,\phi)$ whose coarsening
is an elliptic tail automorphism of order $6$.
\end{teo}

In order to prove this, we start by a lemma about
 an admissible $G'$-cover $F\to C$
over a $2$-marked stable curve $(C;p_1,p_2)$,
where $G'$ is an abelian group.
We observe that any conjugacy class in an abelian
group contains exactly one element, therefore
a $G'$-type (see Definition \ref{def_locgtype}) is an element of $G'$.
Moreover,
 if $\tilde p_i$ is a preimage in $F$ of a
marked point $p_i$, then
the local index at $\tilde p_i$ equals the $G'$-type at $p_i$.

\begin{lemma}\label{lem_2m}
If $G'$ is an abelian group, $(C;p_1,p_2)$ a 
$2$-marked stable curve, and $F\to C$ and admissible $G'$-cover
over $(C;p_1,p_2)$, then the $G'$-types $h_1$ and $h_2$ at $p_1$ and
$p_2$ respectively, are inverses,
$h_1=h_2^{-1}$.
\end{lemma}
\proof We consider at first the case of a smooth $2$-marked curve $(C;p_1,p_2)$.
Because of the monodromy  description given in Proposition \ref{prop_gfg2} and Remark \ref{rmk_mono}, the product $h_1h_2$
is in the commutators subgroup of $G'$, which is trivial because $G'$ is abelian. Therefore
$h_1h_2=1$.

In the case of a general stable curve $C$, 
 we denote by $\tilde p_1^{(i)},\dots,\tilde p_{m_i}^{(i)}$ the marked points on $\overline C_i$, 
\cor{i.e.}~the preimages of $p_1,p_2$ or the $C$ nodes. By the previous point,
if $h_j^{(i)}$ is the $G'$-type of $F$ at the marked point $\tilde p_j^{(i)}$, then
$\prod_{j=1}^{m_i}h_j^{(i)}=1$,
for every $i$. By the balancing condition, for every $G'$-type $h_j^{(i)}$ coming from
a $C$-node, there exists another marked point on $\overline C$ with $G'$-type $h_{j'}^{(i')}=(h_j^{(i)})^{-1}$.
Therefore
$$1=\prod_{j,i}h_j^{(i)}=h_1\cdot h_2.$$\fine

\begin{lemma}
If $(\ssC,\phi)$ is a twisted $S_3$-cover and
 $\msf a\in\underline\Aut_C(\ssC,\phi)\slash\QR(\ssC,\phi)$ is a ghost automorphism,
then $\age(\msf a)\geq 1$.
\end{lemma}

\proof Given a twisted $S_3$-cover $(\ssC,\phi)$, we denote by $F\to C$ the associated
admissible $S_3$-cover and by $\tilde \Gamma\to \Gamma$ the associated graph $S_3$-cover.
We recall that  $b_F$ is the
 index cochain of~$F$.

We prove that if $\msf a$ is a ghost automorphism in $\underline\Aut_C(\ssC,\phi)$
such that $\msf a(e)=1$ for every separating edge of $\Gamma$,
then $\age \msf a\geq 1$. By Lemma \ref{ghost1} this implies the thesis.
By Remark~\ref{rmk_prod},  we have the cycle condition that
for any cycle $(\tilde e_1,\dots,\tilde e_k)$ of $\tilde \Gamma$,
 $\prod (b_F\cdot\msf a)(\tilde e_i)=1$. As any $\msf a(e)$
has order $2$ or $3$ for any $e$, and thus gives an aging of at least $1\slash 2$ or $1\slash 3$ respectively,
 the only case where $\age\msf a<1$ is if there exist two edges
$e_1,e_2\in\E(\Gamma)$ such that $\msf a(e)=1$ if $e\notin\{e_1,e_2\}$ and
$\msf a(e_1)=\msf a(e_2)\in\mmu_3$.
In order to respect the cycles condition, we have a dual graph $\Gamma$ of the type
$$
\xymatrix{
\Gamma'\ar@{-}@/^/[rr]^{e_1}\ar@{-}@/_/[rr]_{e_2}& & \Gamma'',
}
$$
where $\Gamma_1$ and $\Gamma_2$ are two subgraphs of $\Gamma$ such that $\msf a(e)=1$
for every edge in $E(\Gamma_1)$ or $E(\Gamma_2)$. These two subgraphs are
associated to two components $C_1,C_2$ of $C$ such that $C=C_1\cup C_2$ and
they intersect in exactly two nodes $q_1,q_2$, corresponding to edges $e_1,e_2$.

We denote by $\tilde \Gamma_1$ and $\tilde\Gamma_2$ the restrictions
of $\tilde\Gamma$ over $\Gamma_1$ and $\Gamma_2$ respectively.
 If both $\tilde \Gamma_1$ and $\tilde \Gamma_2$ are connected, we denote by
 $\tilde e_1$ and $\tilde e_2$ two preimages of $e_1$ and $e_2$ in $\E(\tilde \Gamma)$
pointing at $\tilde \Gamma_2$ and $\tilde \Gamma_1$ respectively. By the cycle condition,
$(b_F\cdot\msf a)(\tilde e_1)\cdot(b_F\cdot\msf a)(\tilde e_2)=1$, but for the same reason
$(b_F\cdot\msf a)(\tilde e_1)\cdot(b_F\cdot\msf a)(g\cdot\tilde e_2)=1$ for any $g$ in $S_3$,
but this is impossible because $(b_F\cdot\msf a)(\tilde e_2)$ is non-trivial.

If one between $\tilde \Gamma_1$ and $\tilde \Gamma_2$, say the first,
is non-connected, we denote by $\tilde\Gamma'_1,\tilde\Gamma''_1$ its
two components (as $r(e_i)=3$, there are no more than two components).
This means that the restriction $F|_{C'}\to C'$ is an admissible $N$-cover, which means that
$F|_{C'}$ is the union of two admissible $\mmu_3$-covers 
 over the $2$-marked curve $(C_1;p_1,p_2)$.
We denote by $\tilde e_1,\tilde e_2$ the two oriented edges over $e_1$ and $e_2$, both touching $\tilde \Gamma_1'$, and pointing
to $\tilde \Gamma_2$ and $\tilde\Gamma_1'$ respectively.
By Lemma~\ref{lem_2m}, $(b_F\cdot\msf a )(\tilde e_1)=(b_F\cdot\msf a)(\tilde e_2)$ and as $\msf a(e_1)$ has order $3$,
then $(b_F\cdot \msf a)(\tilde e_1)$ and $(b_F\cdot\msf a)(\tilde e_2)$ have order $3$ too.

The oriented edges
$\tilde e_1$ and $\tilde e_2$ touch the same connected components of $\tilde\Gamma''$.
Indeed, if $\tilde \Gamma''$ is non-connected, by local index considerations, both edges
have to touch the same component.
Therefore there exists a cycle passing through $\tilde e_1$ and $\tilde e_2$ and whose other
edges have $\msf a(\tilde e)=1$.

$$
\begin{tikzcd}
\tilde\Gamma'_1\ar[rr, "\tilde e_1", bend left=12]& & \tilde \Gamma''\ar[ll, "\tilde e_2"', bend left=12]\ar[dd, dashed]\ar[lld, bend left=10]\\
\tilde\Gamma'_2\ar[rru, bend left=2] \ar[d,dashed]& &\\
\Gamma'\ar[rr, dash, "e_1", bend left=5] \ar[rr, dash, "e_2"', bend right=5]& &\Gamma''.\\
\end{tikzcd}
$$

Finally, again by the cycle condition we have
$$(b_F\cdot\msf a)(\tilde e_1)\cdot(b_F\cdot\msf a)(\tilde e_2)=(b_F\cdot\msf a)(\tilde e_1)^2=1,$$
but this is a contradiction because $(b_F\cdot\msf a)(e_1)$ has order $3$.\fine

We proved that, as in the case of $G$ abelian group, also for $G=S_3$ the non-canonical locus
$\sing^{\nc}\overline\R_{g,G}$ coincides with the $T$-locus. This is a fundamental result
to approach the extension of pluricanonical forms over a desingularization $\widehat\R_{g,G}\to\overline\R_{g,G}$.

\bibliography{mybiblio}{}
\bibliographystyle{plain} 

\end{document}